\documentclass[11pt,a4paper]{article}

\usepackage{amsmath,amssymb}
\newcommand{\e}{\varepsilon}
\newcommand{\va}{\varphi}
\newcommand{\D}{\Delta}
\newcommand{\La}{\Lambda}
\newcommand{\n}{\nabla}
\newcommand{\N}{\frac{N}{2}}
\newcommand{\NN}{\frac{N}{p}}

\newcommand{\p}{\partial}
\newcommand{\C}{\mathbb{C}}
\newcommand{\R}{\mathbb{R}}
\newcommand{\h}{\hookrightarrow}

\newcommand{\de}{\delta}

\newtheorem{definition}{Definition}
\newtheorem{theorem}{Theorem}
\newtheorem{proposition}{Proposition}
\newtheorem{corollaire}{Corollary}

\newtheorem{notation}{Notation}
\newtheorem{remarka}{Remark}

\newtheorem{lemme}{Lemma}

\title{Existence of global strong solutions in critical spaces for barotropic viscous fluids}
\author{Boris Haspot \thanks{Karls Ruprecht Universit\"at Heidelberg, Institut for Applied Mathematics,
Im Neuenheimer Feld 294,
D-69120 Heildelberg, Germany.
Tel. 49(0)6221-54-6112 }}
\date{}
\begin{document}
%\tableofcontents
\maketitle
%On veut decoupler au niveau de la pression on a alors a etudier:
%$$\p_{t}(u-g)-\frac{1}{\rho}(\mu\D(u-g)+\lambda\n{\rm div}(u-g))=-u\cdot\n u+f+H,$$
%where we set: ${\rm div}g=P(\rho)+h$ such that $\D h=0$ what means that $h$ is very regular and is Fourier transform is support in a compact (in fact we can choose $h=0$).\\
%We have then:
%$$\n{\rm div}g=\n P(\rho)\;\;\;\mbox{and}\;\;\;\D g$$
%En fait on suppose simplement que $\D g=\n P(\rho)$. Dans ce cas on a vraissemblablement: $\n{\rm div}g=\n P(\rho)$.
%En effet on a alors $\D\n{\rm div}g=\D\n P(\rho)$ d'ou notre resultat en fait: $\n{\rm div}g=\n P(\rho)+h$ avec $\D h=0$.
%Si on a des difficultes, on peut loaliser l'operateur.\\
%Cependant ce $h$ eventuellement une constante est genant (au pire on choisit $\lambda=0$ en particulier shallow-water).\\
%\\
%En fait on pose ${\rm div}g=P(\rho)$ alors la solution s'ecrit sous la forme $g_{h}=\n F+h$ avec ${\rm div}h=0$, $F$ est quand a lui unique.
%On a alors $\n {\rm div}g=\n P(\rho)$ et $\D F=P(\rho)$ d'ou $\D\n F=\n P(\rho)=\D g-\D h$.\\
%Il suffit de choisir $g_{0}$ avec $h=0$.
\begin{abstract}
This paper is dedicated to the study of viscous compressible barotropic fluids in dimension $N\geq2$. We address
the question of the global existence of strong solutions for initial data close from a constant state having critical Besov regularity. In a first time, this article show the recent results of 
\cite{CD} and \cite{CMZ} with a new proof. Our result relies on a new a priori estimate for the velocity, where we introduce a new structure to \textit{kill} the coupling between the density and the velocity as in \cite{H2}. We study so a new variable that we call effective velocity. In a second time we improve the results of \cite{CD} and \cite{CMZ} by adding some regularity on the initial data in particular $\rho_{0}$ is in $H^{1}$.  In this case we obtain global strong solutions for a class of large initial data on the density and the velocity which in particular improve the results of D. Hoff in \cite{5H4}. We conclude by generalizing these results for general viscosity coefficients.
\end{abstract}
\section{Introduction}
The motion of a general barotropic compressible fluid is described by the following system:
\begin{equation}
\begin{cases}
\begin{aligned}
&\p_{t}\rho+{\rm div}(\rho u)=0,\\
&\p_{t}(\rho u)+{\rm div}(\rho u\otimes u)-{\rm div}(\mu(\rho)D(u))-\n(\lambda(\rho){\rm div} u)+\n P(\rho)=\rho f,\\
&(\rho,u)_{/t=0}=(\rho_{0},u_{0}).
\end{aligned}
\end{cases}
\label{0.1}
\end{equation}
Here $u=u(t,x)\in\R^{N}$ stands for the velocity field and $\rho=\rho(t,x)\in\R^{+}$ is the density.
The pressure $P$ is a suitable smooth function of $\rho$.
We denote by $\lambda$ and $\mu$ the two viscosity coefficients of the fluid,
which are assumed to satisfy $\mu>0$ and $\lambda+2\mu>0$ (in the sequel to simplify the calculus we will assume the viscosity coefficients as constants). Such a conditions ensures
ellipticity for the momentum equation and is satisfied in the physical cases where $\lambda+\frac{2\mu}{N}>0$.
We supplement the problem with initial condition $(\rho_{0},u_{0})$ and an outer force $f$.
Throughout the paper, we assume that the space variable $x\in\R^{N}$ or to the periodic
box ${\cal T}^{N}_{a}$ with period $a_{i}$, in the i-th direction. We restrict ourselves the case $N\geq2$.\\
The problem of existence of global solution in
time for Navier-Stokes equations was addressed in one dimension for
smooth enough data by Kazhikov and Shelukin in \cite{5K}, and for
discontinuous ones, but still with densities away from zero, by Serre
in  \cite{5S} and Hoff in \cite{5H1}. Those results have been
generalized to higher dimension by Matsumura and Nishida in
\cite{MN} for smooth data close
to equilibrium and by Hoff in the case of discontinuous data in \cite{5H2,5H3}. All
those results do not require to be far from the vacuum.
The existence and uniqueness of local classical solutions for (\ref{0.1})
with smooth initial data such that the density $\rho_{0}$ is bounded
and bounded away from zero (i.e.,
$0<\underline{\rho}\leq\rho_{0}\leq M$)
has been stated by Nash in \cite{Nash}. Let us emphasize that no stability condition was required there.
On the other hand, for small smooth perturbations of a stable
equilibrium with constant positive density, global well-posedness
has been proved in \cite{MN}. Many works on the case of the one dimension have been devoted
to the qualitative behavior of solutions for large time (see for
example \cite{5H1,5K}). Refined functional analysis has been used
for the last decades, ranging from Sobolev, Besov, Lorentz and
Triebel spaces to describe the regularity and long time behavior of
solutions to the compressible model \cite{5So}, \cite{5V},
\cite{5H4}, \cite{5K1}.
Let us recall that (local) existence and uniqueness for (\ref{0.1}) in the case of smooth
data with no vacuum has been stated for long in the pioneering works by J. Nash \cite{Nash},
and A. Matsumura, T. Nishida \cite{MN}.\\
Guided in our approach by numerous works dedicated to the incompressible Navier-Stokes equation (see e.g \cite{Meyer}):
$$
\begin{cases}
\begin{aligned}
&\p_{t}v+v\cdot\n v-\mu\D v+\n\Pi=0,\\
&{\rm div}v=0,
\end{aligned}
\end{cases}
\leqno{(NS)}
$$
we aim at solving (\ref{0.1}) in the case where the data $(\rho_{0},u_{0},f)$ have \textit{critical} regularity.\\
By critical, we mean that we want to solve the system functional spaces with norm
in invariant by the changes of scales which leaves (\ref{0.1}) invariant.
In the case of barotropic fluids, it is easy to see that the transformations:
%$$v(t,x)\longrightarrow lv(l^{2}t,lx),\;\;\;l\in\R$$
%have that property.
%For barotropic fluids, one can check that the appropriate transformation are
\begin{equation}
(\rho(t,x),u(t,x))\longrightarrow (\rho(l^{2}t,lx),lu(l^{2}t,lx)),\;\;\;l\in\R,
\label{1}
\end{equation}
have that property, provided that the pressure term has been changed accordingly.
\begin{definition}
Let $\bar{\rho}>0$. In the sequel we will note: $q=\frac{\rho-\bar{\rho}}{\bar{\rho}}$.
\end{definition}
The use of critical functional frameworks led to several new weel-posedness results for compressible
fluids (see \cite{DL,DW, H1, H2}). In addition to have a norm invariant by (\ref{1}),
appropriate functional space for solving (\ref{0.1}) must provide a control on the $L^{\infty}$
norm of the density (in order to avoid vacuum and loss of ellipticity). For that reason,
we restricted our study to the case where the initial data $(\rho_{0},u_{0})$ and external force $f$
are such that, for some positive constant $\bar{\rho}$:
$$q_{0}\in B^{\NN}_{p,1},\;u_{0}\in B^{\frac{N}{p_{1}}-1}_{p_{1},1}\;\;\mbox{and}\;\;f\in L^{1}_{loc}(\R^{+},\in B^{\frac{N}{p_{1}}-1}_{p_{1},1})$$
with $(p,p_{1})\in [1,+\infty[$ good chosen.\\
Concerning the global existence of strong solutions for initial data with high regularity order and close to a stable equilibrium has been proved by Matsumura and Nishida in \cite{MN} for three-dimensional polytropic ideal fluids and no outer force. More precisely for $\bar{\rho}>0$, the initial data are choosen small in the following spaces $(\rho_{0}-\bar{\rho},u_{0})\in H^{3}\times H^{3}$.
More recently D. Hoff in \cite{5H4,5H5} stated the existence of global weak solutions with small initial data including discontinuous initial data (namely $q_{0}$ is small in $L^{2}\cap L^{\infty}$ and $u_{0}$ is small in $L^{4}$ if $N=2$ and small
in $L^{8}$ if $N=3$). One of the major interest of the results of Hoff is to get some smoothing effects on the incompressible part of the velocity $u$ and on the effective viscous flux $F=(2\mu+\lambda){\rm div}u-P(\rho)+P(\bar{\rho})$ are also pointed out. D. Hoff is the first author to have introduced the notion of effective flux which play a crucial role in the proof of P-L Lions for the existence of global weak solution. However if the results of Hoff are critical in the sense of the scaling for the density, it is not the case for the initial velocity. In \cite{5H2}, D. Hoff show a very interesting theorem of weak-strong uniqueness when $P(\rho)=K\rho$ with $K>0$.
To speak roughly under the conditions that two solutions $(\rho,u)$, $(\rho_{1},u_{1})$ check a control $L^{\infty}$ on the density and a control Lipschitz on the velocity, with additional property of regularity on the strong solution $(\rho_{1},u_{1})$ then we obtain $(\rho,u)=(\rho_{1},u_{1})$.
D. Hoff use this result to show that the solutions of \cite{5H2} are unique. We will use this theorem in the sequel, by showing that our solutions verify the hypothesis of D. Hoff in \cite{5H2}.
Finally R. Danchin in \cite{DG} show for the first time a result of existence of global strong solution close from a stable equilibrium in critical space for the scaling of the system. More precisely the initial data are choose as follows $(q_{0},u_{0})\in (B^{\N}_{2,1}\cap B^{\N-1}_{2,1})\times B^{\N-1}_{2,1}$. The main difficulty is to get estimates on the linearized system where the velocity and the density are coupled via the pressure, and what is crucial in this work is the smoothing effect on the velocity and a $L^{1}$ decay on $\rho-\bar{\rho}$ (this play a necessary role to control the pressure term). In this work, R. Danchin use some astucious inequality of energy on the system in variable Fourier where he has decomposed the space in dyadic shell. This explain in particular why the result is obtained in Besov space with a Lebesgue index $p=2$. In the same time of the redaction of this paper, Q. Chen et al in \cite{CMZ} and F. Charve and R. Danchin in \cite{CD} improve the previous result by working in more general Besov space by studying the linear part of the system.\\
The goal of this article is to make a connection between the article of D. Hoff \cite{5H4,5H5}
and those of Q. Chen et al and F. Charve and R. Danchin in \cite{CD} and \cite{CMZ}. In fact we extend the results \cite{CD} and \cite{CMZ} to the case where the Lebesgue index of Besov space are more general, it means $q_{0}\in B^{\NN}_{p,1}$ and $\;u_{0}\in B^{\frac{N}{p_{1}}-1}_{p_{1},1}$ with $p$ and $p_{1}$ good choosen. In \cite{CD} and \cite{CMZ}, the authors  obtain global weak solutions  when $p=p_{1}<2N$ and strong solution when $p=p_{1}\leq N$, the restriction on the choice of $p$ come from a very strong coupling between the pressure
and the velocity. Indeed in this case the coupling is very strong in high frequencies because of the term of pressure, that's why we need to integrate completly the pressure term in the linear part. In the case of lows frequencies, according the point of view of the Fourier frequencies, the term of pressure is very regular so it does not make problem to consider in the rest. The study of the linear part in this case is crucial to get a gain $L^{1}$ of integrability on the density, and for low frequencies we follow the method of R. Danchin in \cite{DG}. In \cite{5H4,5H5}, D. Hoff get global weak solution with a critical regularity on the density in the sense that $\rho_{0}-\bar{\rho}$ is small in $L^{\infty}$. It means that he does not ask any regularity on the initial density which is of this point of view besser than \cite{DG}, \cite{CD} and \cite{CMZ}. However the velocity in his case in only $L^{1}$ $log$ Lipschitz, that is why he can not obtain the uniqueness and have only weak solution. In this paper we improve the results of \cite{5H4,5H5} by the fact that we get strong solutions, and we will show that it is just enough to ask a arbitrary small $\e$ regularity on the initial density to get uniqueness.\\
%To be more precise, the pressure term is considered as a term of rest for the elliptic operator in the momentum equation of (\ref{0.1}).
In \cite{H2}, we improve the results of R. Danchin in \cite{DL,DW}, in the sense that the initial density belongs to larger spaces $B^{\NN}_{p,1}$ with $p\in[1,+\infty[$. In the present paper, we address the question of global existence of strong solution in the critical functional framework under the assumption
that the initial density belongs to critical Besov space with a index of integrability different of this of the velocity. To do that, as in \cite{H2} we  introduce a new variable in high frequencies than the velocity that we call effective velocity in the goal to \textit{kill}
the relation of coupling between the velocity and the pressure. We observe that this new notion of effective velocity allow us easily to get as R. Danchin in \cite{DG} a $L^{1}$ decay on $q$. However this new variable is interesting only in high frequencies, indeed in low frequencies the term $\n P(\rho)$ is small in Fourier analysis. Moreover in the low frequency regime, the first order terms predominate and the viscous term $\D u$ may be neglected in Fourier analysis, so that (\ref{0.1}) has to be treat by means
of hyperbolic energy methods (more particularly the velocity verifies in some way a wave equation). This implies that we can treat the low regime only in space construct on $L^{2}$, it is classical that the hyperbolic system are ill-posed in general $L^{p}$ spaces.
%to make the additional assumption that $\rho-\bar{\rho}$ is {\it small} in $B^{\NN}_{p,1}$.
%The reason why is that we handled the elliptic operator in the momentum equation of $(SW)$
%as a constant coefficient second order operator plus a perturbation introduced by $\rho-\bar{\rho}$ which, if sufficiently
%small, may be treated as a harmless source term. For smoother data however, the additional regularity can compensate {\it large}
%perturbations of the constant coefficient operator. This fact has been used in \cite{DL} leads to local well-posedness
%results for smooth enough data with density bounded away from zero. The price to pay however is that assuming extra smoothness
%precludes from using a critical framework.\\
So as in \cite{CMZ} and \cite{CD}, the system has to be handled differently in low and high frequencies. In short, we will use the analysis of R. Danchin in \cite{DG} in low frequencies and the introduction of this new variable the effctive velocity introduced in \cite{H2} in high frequencies.
%We adapt the spirit of the results %of \cite{AP} and
%of \cite{H2} which treat the case of strong solution in finite time.
To simplify the notation, we assume from now on that $\bar{\rho}=1$. Hence as long as $\rho$ does not vanish, the equations for ($q=\rho-1$,$u$) read:
\begin{equation}
\begin{cases}
\begin{aligned}
&\p_{t}q+u\cdot\n q=-(1+q){\rm div}u,\\
&\p_{t}u+u\cdot\n u-\frac{1}{1+q}{\cal A}u+\n P(1+q)=f,
\end{aligned}
\end{cases}
\label{0.6}
\end{equation}
In the sequel we will note ${\cal A}=\mu\D+(\lambda+\mu)\n{\rm div}$ and where $g$ is a smooth function which may be computed from the pressure function $P$.
One can now state our main result.
\begin{theorem}
Let $P$ a suitably smooth function of the density such that $P^{'}(1)>0$, $f\in\widetilde{L}^{1}(\widetilde{B}^{\N-1,\frac{N}{p_{1}}-1}_{2,p_{1},1})$ and $1\leq p_{1}\leq p<+\infty$
such that $\frac{1}{p_{1}}\leq\frac{1}{N}+\frac{1}{p}$.
%and $2\NN-1>0$.
Assume that $u_{0}\in \widetilde{B}^{\N-1,\frac{N}{p_{1}}-1}_{2,p_{1},1}$, $f\in L^{1}_{loc}(\R^{+},\widetilde{B}^{\N-1,\frac{N}{p_{1}}-1}_{2,p_{1},1})$ and $q_{0}\in \widetilde{B}^{\N-1,\NN}_{2,p,1}$.%\cap B^{\frac{N}{p_{1}}-2}_{p_{1},1}$,
Then there exists a constant $\e_{0}$ such that if:
$$\|q_{0}\|_{\widetilde{B}^{\N-1,\NN}_{2,p,1}}+\|u_{0}\|_{\widetilde{B}^{\N-1,\frac{N}{p_{1}}-1}_{2,p_{1},1}}+\|f\|_{\widetilde{L}^{1}(\widetilde{B}^{\N-1,\frac{N}{p_{1}}-1}_{2,p_{1},1})}\leq\e_{0},$$
then if $\frac{1}{p}+\frac{1}{p_{1}}>\frac{1}{N}$, $p<\max(4,N)$ and $\frac{1}{2}\leq\frac{1}{p}+\frac{1}{p_{1}}$ there exists a global solution $(q,u)$ for system (\ref{0.1})
with $1+q$ bounded away from zero and,
$$
\begin{aligned}
&q\in \widetilde{C}(\R,\widetilde{B}^{\N-1,\NN}_{2,p,1}%\cap B^{\frac{N}{p_{2}}-1}_{p_{1},1}
)\cap\widetilde{L}^{1}(\R,
\widetilde{B}^{\N+1,\frac{N}{p}}_{2,p,1})\;\;\;\mbox{and}\\
&\hspace{2cm}\;u\in \widetilde{C}(\R;\widetilde{B}^{\N-1,\frac{N}{p_{1}}-1}_{2,p_{1},1}+\widetilde{B}^{\N-1,\NN}_{2,p,1})
\cap\widetilde{L}^{1}(\R,\widetilde{B}^{\N+1,\NN+1}_{2,p,1}).
\end{aligned}
$$
Moreover this solution is unique if $\frac{2}{N}\leq\frac{1}{p}+\frac{1}{p_{1}}$.
\label{theo1}
\end{theorem}
\begin{remarka}
 This theorem is the same than Chen et al obtained in \cite{CMZ}. In particular we have strong restrictions on $p$, it means $p<\max(4,2N)$. This fact is due to the interactions between low and high frequencies in the paraproduct laws. However we obtain this result with a new method which seems more flexible and we will explain why in thecorrolary \ref{corollaire1}.
\end{remarka}

\begin{remarka}
It seems possible to improve the theorem \ref{theo1} by choosing initial data $q_{0}$ in $B^{\N-1,\NN}_{(2,1),(p,\infty)}\cap B^{\N-1,0}_{(2,1),(\infty,1)}$, however some supplementary conditions appear on $p_{1}$ in this case. Here $B^{s_{1},s_{2}}_{(p_{1},r_{1}),(p_{2},r_{2})}$ is a Besov space where the behavior is $B^{s_{1}}_{p_{1},r_{1}}$  and $B^{s_{2}}_{p_{2},r_{2}}$  in high frequencies.
\end{remarka}
The key to theorem \ref{theo1} is to introduce a new variable $v_{1}$ to control the velocity where to avoid the coupling between the density and the velocity, we analyze by a new way the pressure term.  More precisely we write the gradient of the pressure as a Laplacian of the variable $v_{1}$, and we introduce this term in the linear part of the momentum equation. We have then a control on $v_{1}$ which can write roughly as $u-{\cal G}P(\rho)$ where ${\cal G}$ is a pseudodifferential operator of order $-1$. We will call $u-{\cal G}P(\rho)$ the effective velocity. By this way, we have canceled the coupling between $v_{1}$ and the density, we next verify easily that we have a control Lipschitz of the gradient of $u$ (it is crucial to estimate the density by the mass equation).
%a new estimate for a class of parabolic systems with coefficients in
%$C([0,T];B^{\NN}_{p,1})$ which is obtained when linearizing the momentum equation in $(SW)$ (see proposition
%\ref{linearise}).
%The basic idea is
%It turns out that our study of the linearization of (\ref{0.1}) leads to also to the following continuation criterion:
%We can now state pur second main result.
In the previous theorem \ref{theo1}, we have as in \cite{CMZ} very big restrictions on $p$ ($p<\max(4,2N)$) because the behavior in low frequencies. At the difference with the results of strong solutions in finite time (see \cite{H2}), we can not choose $p$ arbitrarly big.
To overcome this difficulty, we need to add some additional conditions on $(q_{0},u_{0})$ in low frequencies as in \cite{CD} to avoid these restrictions in the use of the paraproduct laws. We obtain then the following corollary:
\begin{corollaire}
Let $P$ a suitably smooth function of the density with $P^{'}(1)>0$, $f\in\widetilde{L}^{1}(\widetilde{B}^{\N-1,\frac{N}{p_{1}}-1}_{2,p_{1},1}\cap B^{0}_{2,r})$ and $1\leq p_{1}\leq p<+\infty$
such that $\frac{1}{p_{1}}\leq\frac{1}{N}+\frac{1}{p}$.
%and $2\NN-1>0$.
Assume that $u_{0}\in \widetilde{B}^{\N-1,\frac{N}{p_{1}}-1}_{2,p_{1},1}\cap B^{0}_{2,r}$, $f\in L^{1}_{loc}(\R^{+},\widetilde{B}^{\N-1,\frac{N}{p_{1}}-1}_{2,p_{1},1})$ and $q_{0}\in \widetilde{B}^{\N-1,\NN}_{2,p,1}\cap B^{0,1}_{2,r}$ with $r=+\infty$ if $N\geq 3$ and $r=1$ if $N=2$.%\cap B^{\frac{N}{p_{1}}-2}_{p_{1},1}$,
Then there exists a constant $\e_{0}$ such that if:
$$\|q_{0}\|_{\widetilde{B}^{\N-1,\NN}_{2,p,1}\cap B^{0,1}_{2,r}}+\|u_{0}\|_{\widetilde{B}^{\N-1,\frac{N}{p_{1}}-1}_{2,p_{1},1}\cap B^{0}_{2,r}}+\|f\|_{\widetilde{L}^{1}(\widetilde{B}^{\N-1,\frac{N}{p_{1}}-1}_{2,p_{1},1}\cap B^{0}_{2,r})}\leq\e_{0},$$
then if $\frac{1}{p}+\frac{1}{p_{1}}>\frac{1}{N}$, there exists a global solution $(q,u)$ for system (\ref{0.1})
with $1+q$ bounded away from zero and,
$$
\begin{aligned}
&q\in \widetilde{C}(\R,\widetilde{B}^{\N-1,\NN}_{2,p,1}\cap B^{0,1}_{2,r}%\cap B^{\frac{N}{p_{2}}-1}_{p_{1},1}
)\cap\widetilde{L}^{1}(\R,
\widetilde{B}^{\N+1,\frac{N}{p}}_{2,p,1}\cap B^{2,1}_{2,r})\;\;\;\mbox{and}\\
&\hspace{2cm}\;u\in \widetilde{C}\big(\R;(\widetilde{B}^{\N-1,\frac{N}{p_{1}}-1}_{2,p_{1},1}+\widetilde{B}^{\N-1,\NN}_{2,p,1})\cap B^{0}_{2,r}\big)
\cap\widetilde{L}^{1}(\R,\widetilde{B}^{\N+1,\NN+1}_{2,p,1}\cap B^{2}_{2,r}).
\end{aligned}
$$
Moreover this solution is unique if $\frac{2}{N}\leq\frac{1}{p}+\frac{1}{p_{1}}$.
\label{corollaire1}
\end{corollaire}
\begin{remarka}
We can observe that when $p$ tends to infinity, we are close to get weak solution with the following initial data  $(q_{0},u_{0})$ in
$B^{\N-1,0}_{2,\infty,1}\times B^{\N-1,1}_{2,N,1}$. It means that this theorem rely the result of D. Hoff where the initial density is assumed $L^{\infty}$ but where the initial velocity is more regular (it means not critical) and the results of R. Danchin in \cite{DW}.Moreover it is right for general pressure when in the works of D. Hoff the pressure verify $P(\rho)=K\rho$ with $K>0$.
\end{remarka}
\begin{remarka}
If $r=+\infty$, then we replace above the strong continuity in $B^{s}_{2,r}$ by the weak continuity.
\end{remarka}
\begin{remarka}
In some some sense, we could consider that the case $p>N$ is not so important because we obtain only the existence of global weak solution as in the works of D. Hoff in \cite{5H4}, \cite{5H2}. However it stays very interesting, indeed as in the work of F. Charve and R. Danchin in \cite{CD}, we could add some additional condition on the data such $u_{0}\in B^{0}_{N,1}$. In this fact, it would be easy to show some results of persistency as for Navier-Stokes without condition of smallness on $\|u_{0}\|_{B^{0}_{N,1}}$, in particular the fact that $u\in\widetilde{L}^{\infty}(\widetilde{B}^{\N-1,0}_{2,N,1}+\widetilde{B}^{\N-1,\NN}_{2,p,1})\cap \widetilde{L}^{1}(\widetilde{B}^{\N+1,\NN+1}_{2,p,1})$. And in this case, we will case the existence of global strong solutions
but with only a condition of smallness for the density on $\|q_{0}\|_{\widetilde{B}^{\N-1,\N}_{2,p,1}\cap \widetilde{B}^{0,1}_{2,r}}$.
In particular it improves widely the results of \cite{CMZ} and \cite{CD}, this result is close from the results of D. Hoff in \cite{5H4} and \cite{5H2} except that we are in critical space for the scaling (except concerning the fact that $q_{0}\in \widetilde{B}^{0,1}_{2,r}$).\\
In particular as in \cite{CD}, we can take $u_{0}(x)=\phi(x)sin(\e^{-1}x\cdot\omega)n$ where $\omega$ and $n$ stand for any unit vectors of $\R^{N}$ and $\phi$ for any smooth compactly supported function then we have if $p_{1}>N$:
$$\|u_{0}\|_{B^{\frac{N}{p_{1}}-1}_{p_{1},1}}\leq C\e^{1-\frac{N}{p_{1}}},$$
so that the smallness condition is satisfied by $u_{0}$ if $\e$ small enough. However we remark that $u_{0}$ is arbitrarly big in $L^{3}$. On the other hand, $u_{0}$ belongs to Scwartz class ${\cal S}$ hence also to $B^{0}_{N,1}$ so that uniqueness holds true by persistency results.
\end{remarka}
\begin{remarka}
 We can observe that $u_{0}\in B^{0}_{2,r}$ corresponds exactly to the energy space when $r=2$, in this sense this additional regularity on the velocity seems very natural. We will explain in the corollary \ref{corollaire2} why it is perfectly adapted in the case of specific viscosity coefficients. Indeed in the general case $q_{0}\in B^{1}_{2,r}$ is not in the energy space.
\end{remarka}
We want treat now the special case of the \textit{BD viscosity coefficients}. Indeed in \cite{BD} Bresch and Desjardins have discovered a new entropy inequality when in (\ref{0.1}), we have:
$$\lambda(\rho)=\rho\mu^{'}(\rho)-\mu(\rho).$$
In this case they show that we can control $\sqrt{\rho}\n\va(\rho)$ in $L^{\infty}(L^{2})$ where $\va^{'}(\rho)=\frac{\mu^{'}(\rho)}{\rho}$. Roughly it means that we controll the density $\rho$ in $L^{\infty}(H^{1})$.
It is the additional condition that we ask in the corollary \ref{corollaire1}. In the following result, we prove that we can extend
the corollary \ref{corollaire1} to the case of general viscosity.
\begin{corollaire}
Let $P$ a suitably smooth function of the density with $P^{'}(1)>0$, $f\in\widetilde{L}^{1}(\widetilde{B}^{\N-1,\frac{N}{p_{1}}-1}_{2,p_{1},1}\cap B^{0}_{2,r})$ and $\mu$, $\lambda$ are general regular functions such that $\mu(1)>0$ and $\mu(1)+\lambda(1)>0$ and $1\leq p_{1}\leq p<+\infty$
such that $\frac{1}{p_{1}}\leq\frac{1}{N}+\frac{1}{p}$.
%and $2\NN-1>0$.
Assume that $u_{0}\in \widetilde{B}^{\N-1,\frac{N}{p_{1}}-1}_{2,p_{1},1}\cap B^{0}_{2,r}$, $f\in L^{1}_{loc}(\R^{+},\widetilde{B}^{\N-1,\frac{N}{p_{1}}-1}_{2,p_{1},1})$ and $q_{0}\in \widetilde{B}^{\N-1,\NN}_{2,p,1}\cap B^{0,1}_{2,r}$ with $r=+\infty$ if $N\geq 3$ and $r=1$ if $N=2$.%\cap B^{\frac{N}{p_{1}}-2}_{p_{1},1}$,
Then there exists a constant $\e_{0}$ such that if:
$$\|q_{0}\|_{\widetilde{B}^{\N-1,\NN}_{2,p,1}\cap B^{0,1}_{2,r}}+\|u_{0}\|_{\widetilde{B}^{\N-1,\frac{N}{p_{1}}-1}_{2,p_{1},1}\cap B^{0}_{2,r}}+\|f\|_{\widetilde{L}^{1}(\widetilde{B}^{\N-1,\frac{N}{p_{1}}-1}_{2,p_{1},1}\cap B^{0}_{2,r})}\leq\e_{0},$$
then if $\frac{1}{p}+\frac{1}{p_{1}}>\frac{1}{N}$, there exists a global solution $(q,u)$ for system (\ref{0.1})
with $1+q$ bounded away from zero and,
$$
\begin{aligned}
&q\in \widetilde{C}(\R,\widetilde{B}^{\N-1,\NN}_{2,p,1}\cap B^{0,1}_{2,r}%\cap B^{\frac{N}{p_{2}}-1}_{p_{1},1}
)\cap\widetilde{L}^{1}(\R,
\widetilde{B}^{\N+1,\frac{N}{p}}_{2,p,1}\cap B^{2,1}_{2,r})\;\;\;\mbox{and}\\
&\hspace{2cm}\;u\in \widetilde{C}\big(\R;(\widetilde{B}^{\N-1,\frac{N}{p_{1}}-1}_{2,p_{1},1}+\widetilde{B}^{\N-1,\NN}_{2,p,1})\cap B^{0}_{2,r}\big)
\cap\widetilde{L}^{1}(\R,\widetilde{B}^{\N+1,\NN+1}_{2,p,1}\cap B^{2}_{2,r}).
\end{aligned}
$$
Moreover this solution is unique if $\frac{2}{N}\leq\frac{1}{p}+\frac{1}{p_{1}}$.
\label{corollaire2}
\end{corollaire}
\begin{remarka}
 This result is very interesting in the case of the \textit{BD viscosity coefficients}. In this case our result is very close of the energy initial data with the optimal condition for the scaling $(q_{0},u_{0})\in B^{0}_{\infty,1}\times B^{0}_{N,1})$. In particular it concerns the shallow-water system.
\end{remarka}
\begin{remarka}
 Moreover our method is more flexible than the proofs of D. Hoff in \cite{5H2}, \cite {5H3}, \cite{5H4} as these works  are based crucially on the notion of effective pressure and on a gain of integrability on the velocity which are right only in the case of constant viscosity coefficients.
\end{remarka}
Our paper is structured as follows. In section \ref{section2}, we give a few notation and briefly introduce the basic Fourier analysis
techniques needed to prove our result. In section \ref{section3}, we prove estimate on the transport equation. In section \ref{section4}, we prove the theorem \ref{theo1}.In section \ref{section5} we prove the corollaries \ref{corollaire1} and \ref{corollaire2}. Two inescapable technical commutator estimates  and the proof of paraproduct in hybrid Besov spaces%and some theorems of ellipticity
are postponed in appendix.
\section{Littlewood-Paley theory and Besov spaces}
\label{section2}
Throughout the paper, $C$ stands for a constant whose exact meaning depends on the context. The notation $A\lesssim B$ means
that $A\leq CB$.
For all Banach space $X$, we denote by $C([0,T],X)$ the set of continuous functions on $[0,T]$ with values in $X$.
For $p\in[1,+\infty]$, the notation $L^{p}(0,T,X)$ or $L^{p}_{T}(X)$ stands for the set of measurable functions on $(0,T)$
with values in $X$ such that $t\rightarrow\|f(t)\|_{X}$ belongs to $L^{p}(0,T)$.
Littlewood-Paley decomposition  corresponds to a dyadic
decomposition  of the space in Fourier variables.
We can use for instance any $\varphi\in C^{\infty}(\R^{N})$,
supported in
${\cal{C}}=\{\xi\in\R^{N}/\frac{3}{4}\leq|\xi|\leq\frac{8}{3}\}$
such that:
$$\sum_{l\in\mathbb{Z}}\varphi(2^{-l}\xi)=1\,\,\,\,\mbox{if}\,\,\,\,\xi\ne 0.$$
Denoting $h={\cal{F}}^{-1}\varphi$, we then define the dyadic
blocks by:
$$\D_{l}u=\varphi(2^{-l}D)u=2^{lN}\int_{\R^{N}}h(2^{l}y)u(x-y)dy\,\,\,\,\mbox{and}\,\,\,S_{l}u=\sum_{k\leq
l-1}\D_{k}u\,.$$ Formally, one can write that:
$$u=\sum_{k\in\mathbb{Z}}\D_{k}u\,.$$
This decomposition is called homogeneous Littlewood-Paley
decomposition. Let us observe that the above formal equality does
not hold in ${\cal{S}}^{'}(\R^{N})$ for two reasons:
\begin{enumerate}
\item The right hand-side does not necessarily converge in
${\cal{S}}^{'}(\R^{N})$.
\item Even if it does, the equality is not
always true in ${\cal{S}}^{'}(\R^{N})$ (consider the case of the
polynomials).
\end{enumerate}
\subsection{Homogeneous Besov spaces and first properties}
\begin{definition}
For
$s\in\R,\,\,p\in[1,+\infty],\,\,q\in[1,+\infty],\,\,\mbox{and}\,\,u\in{\cal{S}}^{'}(\R^{N})$
we set:
$$\|u\|_{B^{s}_{p,q}}=(\sum_{l\in\mathbb{Z}}(2^{ls}\|\D_{l}u\|_{L^{p}})^{q})^{\frac{1}{q}}.$$
The Besov space $B^{s}_{p,q}$ is the set of temperate distribution $u$ such that $\|u\|_{B^{s}_{p,q}}<+\infty$.
\end{definition}
%A difficulty due to the choice of homogeneous spaces arises at this
%point. Indeed, $ \|.\|_{B^{s}_{p,q}}$ cannot be a norm on
%$\{u\in{\cal{S}}^{'}(\R^{N}),\|u\|_{B^{s}_{p,q}}<+\infty\}$ because
%$\|u\|_{B^{s}_{p,q}}=0$ means that $u$ is a polynomial. This
%enforces us to adopt the following definition for homogeneous Besov
%spaces, see \cite{Bou}.
%\begin{definition}
%Let $s\in\R,\,\,p\in[1,+\infty],\,\,q\in[1,+\infty]$.\\
%Denote $m=[s-N/p]$ if $s-N/p\notin\mathbb{Z}$ or $q>1$ and
%$m=s-N/p-1$ otherwise.
%\begin{itemize}
%\item If $m<0$, then we define $B^{s}_{p,q}$ as:
%$$B^{s}_{p,q}=\biggl\{u\in{\cal{S}}^{'}(\R^{N})\;\;/\;\;\|u\|_{B^{s}_{p,q}}<\infty\,\,\,\mbox{and}
%\,\,\,u=\sum_{l\in\mathbb{Z}}\D_{l}u\,\,\mbox{in}\,\,{\cal{S}}^{'}(\R^{N})\biggl\}\,.$$
%\item If $m\geq 0$, we denote by ${\cal{P}}_{m}[\R^{N}]$ the set of
%polynomials of degree less than or equal to $m$ and we set:
%$$B^{s}_{p,q}=\biggl\{u\in{\cal{S}}^{'}(\R^{N})/{\cal{P}}_{m}[\R^{N}]\;\;/\;\;
%\|u\|_{B^{s}_{p,q}}<\infty\,\,\,and\,\,\,u=\sum_{l\in\mathbb{Z}}
%\D_{l}u\,\,in\,\,{\cal{S}}^{'}(\R^{N}){\cal{P}}_{m}[\R^{N}]\biggl\}\,.$$
%\end{itemize}
%\end{definition}
\begin{remarka}The above definition is a natural generalization of the
nonhomogeneous Sobolev and H$\ddot{\mbox{o}}$lder spaces: one can show
that $B^{s}_{\infty,\infty}$ is the nonhomogeneous
H$\ddot{\mbox{o}}$lder space $C^{s}$ and that $B^{s}_{2,2}$ is
the nonhomogeneous space $H^{s}$.
\end{remarka}
\begin{proposition}
\label{derivation,interpolation}
The following properties holds:
\begin{enumerate}
\item there exists a constant universal $C$
such that:\\
$C^{-1}\|u\|_{B^{s}_{p,r}}\leq\|\n u\|_{B^{s-1}_{p,r}}\leq
C\|u\|_{B^{s}_{p,r}}.$
\item If
$p_{1}<p_{2}$ and $r_{1}\leq r_{2}$ then $B^{s}_{p_{1},r_{1}}\hookrightarrow
B^{s-N(1/p_{1}-1/p_{2})}_{p_{2},r_{2}}$.
\item $B^{s^{'}}_{p,r_{1}}\hookrightarrow B^{s}_{p,r}$ if $s^{'}> s$ or if $s=s^{'}$ and $r_{1}\leq r$.
 %$(B^{s_{1}}_{p,r},B^{s_{2}}_{p,r})_{\theta,r^{'}}=B^{\theta
%s_{1}+(1-\theta)s_{2}}_{p,r^{'}}$.
\end{enumerate}
\label{interpolation}
\end{proposition}
%Before going further into the paraproduct for Besov spaces, let us state an important proposition.
%\begin{proposition}
%Let $s\in\R$ and $1\leq p,r\leq+\infty$. Let $(u_{q})_{q\geq-1}$ be a sequence of functions such that
%$$(\sum_{q\geq-1}2^{qsr}\|u_{q}\|_{L^{p}}^{r})^{\frac{1}{r}}<+\infty.$$
%If $\mbox{supp}\hat{u}_{1}\subset {\cal C}(0,2^{q}R_{1},2^{q}R_{2})$ for some $0<R_{1}<R_{2}$ then $u=\sum_{q\geq-1}u_{q}$ belongs to $B^{s}_{p,r}$ and there exists a universal constant $C$ such that:
%$$\|u\|_{B^{s}_{p,r}}\leq C^{1+|s|}\big(\sum_{q\geq-1}(2^{qs}\|u_{q}\|_{L^{p}})^{r}\big)^{\frac{1}{r}}.$$
%\label{resteimp1}
%\end{proposition}
Let now recall a few product laws in Besov spaces coming directly from the paradifferential calculus of J-M. Bony
(see \cite{BJM}) and rewrite on a generalized form in \cite{AP} by H. Abidi and M. Paicu (in this article the results are written
in the case of homogeneous sapces but it can easily generalize for the nonhomogeneous Besov spaces).
\begin{proposition}
\label{produit1}
We have the following laws of product:
\begin{itemize}
\item For all $s\in\R$, $(p,r)\in[1,+\infty]^{2}$ we have:
\begin{equation}
\|uv\|_{B^{s}_{p,r}}\leq
C(\|u\|_{L^{\infty}}\|v\|_{B^{s}_{p,r}}+\|v\|_{L^{\infty}}\|u\|_{B^{s}_{p,r}})\,.
\label{2.2}
\end{equation}
\item Let $(p,p_{1},p_{2},r,\lambda_{1},\lambda_{2})\in[1,+\infty]^{2}$ such that:$\frac{1}{p}\leq\frac{1}{p_{1}}+\frac{1}{p_{2}}$,
$p_{1}\leq\lambda_{2}$, $p_{2}\leq\lambda_{1}$, $\frac{1}{p}\leq\frac{1}{p_{1}}+\frac{1}{\lambda_{1}}$ and
$\frac{1}{p}\leq\frac{1}{p_{2}}+\frac{1}{\lambda_{2}}$. We have then the following inequalities:\\
if $s_{1}+s_{2}+N\inf(0,1-\frac{1}{p_{1}}-\frac{1}{p_{2}})>0$, $s_{1}+\frac{N}{\lambda_{2}}<\frac{N}{p_{1}}$ and
$s_{2}+\frac{N}{\lambda_{1}}<\frac{N}{p_{2}}$ then:
\begin{equation}
\|uv\|_{B^{s_{1}+s_{2}-N(\frac{1}{p_{1}}+\frac{1}{p_{2}}-\frac{1}{p})}_{p,r}}\lesssim\|u\|_{B^{s_{1}}_{p_{1},r}}
\|v\|_{B^{s_{2}}_{p_{2},\infty}},
\label{2.3}
\end{equation}
when $s_{1}+\frac{N}{\lambda_{2}}=\frac{N}{p_{1}}$ (resp $s_{2}+\frac{N}{\lambda_{1}}=\frac{N}{p_{2}}$) we replace
$\|u\|_{B^{s_{1}}_{p_{1},r}}\|v\|_{B^{s_{2}}_{p_{2},\infty}}$ (resp $\|v\|_{B^{s_{2}}_{p_{2},\infty}}$) by
$\|u\|_{B^{s_{1}}_{p_{1},1}}\|v\|_{B^{s_{2}}_{p_{2},r}}$ (resp $\|v\|_{B^{s_{2}}_{p_{2},\infty}\cap L^{\infty}}$),
if $s_{1}+\frac{N}{\lambda_{2}}=\frac{N}{p_{1}}$ and $s_{2}+\frac{N}{\lambda_{1}}=\frac{N}{p_{2}}$ we take $r=1$.
\\
If $s_{1}+s_{2}=0$, $s_{1}\in(\frac{N}{\lambda_{1}}-\frac{N}{p_{2}},\frac{N}{p_{1}}-\frac{N}{\lambda_{2}}]$ and
$\frac{1}{p_{1}}+\frac{1}{p_{2}}\leq 1$ then:
\begin{equation}
\|uv\|_{B^{-N(\frac{1}{p_{1}}+\frac{1}{p_{2}}-\frac{1}{p})}_{p,\infty}}\lesssim\|u\|_{B^{s_{1}}_{p_{1},1}}
\|v\|_{B^{s_{2}}_{p_{2},\infty}}.
\label{2.4}
\end{equation}
If $|s|<\NN$ for $p\geq2$ and $-\frac{N}{p^{'}}<s<\NN$ else, we have:
\begin{equation}
\|uv\|_{B^{s}_{p,r}}\leq C\|u\|_{B^{s}_{p,r}}\|v\|_{B^{\NN}_{p,\infty}\cap L^{\infty}}.
\label{2.5}
\end{equation}
\end{itemize}
\end{proposition}
\begin{remarka}
In the sequel $p$ will be either $p_{1}$ or $p_{2}$ and in this case $\frac{1}{\lambda}=\frac{1}{p_{1}}-\frac{1}{p_{2}}$
if $p_{1}\leq p_{2}$, resp $\frac{1}{\lambda}=\frac{1}{p_{2}}-\frac{1}{p_{1}}$
if $p_{2}\leq p_{1}$.
\end{remarka}
\begin{corollaire}
\label{produit2}
Let $r\in [1,+\infty]$, $1\leq p\leq p_{1}\leq +\infty$ and $s$ such that:
\begin{itemize}
\item $s\in(-\frac{N}{p_{1}},\frac{N}{p_{1}})$ if $\frac{1}{p}+\frac{1}{p_{1}}\leq 1$,
\item $s\in(-\frac{N}{p_{1}}+N(\frac{1}{p}+\frac{1}{p_{1}}-1),\frac{N}{p_{1}})$ if $\frac{1}{p}+\frac{1}{p_{1}}> 1$,
\end{itemize}
then we have if $u\in B^{s}_{p,r}$ and $v\in B^{\frac{N}{p_{1}}}_{p_{1},\infty}\cap L^{\infty}$:
$$\|uv\|_{B^{s}_{p,r}}\leq C\|u\|_{B^{s}_{p,r}}\|v\|_{B^{\frac{N}{p_{1}}}_{p_{1},\infty}\cap L^{\infty}}.$$
\end{corollaire}
%For a proof of this proposition see \cite{5DG}. The limit case $s_{1}+s_{2}=t_{1}+t_{2}=0$ in (\ref{5prodinteressant2}) is of interest.
%When $p\geq2$, the following estimate holds true whenever $s$ is in the range $(-\NN,\NN]$ (see e.g. \cite{5RS}):
%\begin{equation}
%\|uv\|_{B^{-\NN}_{p,\infty}}\leq C\|u\|_{B^{s}_{p,1}}\|v\|_{B^{-s}_{p,\infty}}.
%\end{equation}
The study of non stationary PDE's requires space of type $L^{\rho}(0,T,X)$ for appropriate Banach spaces $X$. In our case, we
expect $X$ to be a Besov space, so that it is natural to localize the equation through Littlewood-Payley decomposition. But, in doing so, we obtain
bounds in spaces which are not type $L^{\rho}(0,T,X)$ (except if $r=p$).
We are now going to
define the spaces of Chemin-Lerner in which we will work, which are
a refinement of the spaces
$L_{T}^{\rho}(B^{s}_{p,r})$.
$\hspace{15cm}$
\begin{definition}
Let $\rho\in[1,+\infty]$, $T\in[1,+\infty]$ and $s_{1}\in\R$. We set:
$$\|u\|_{\widetilde{L}^{\rho}_{T}(B^{s_{1}}_{p,r})}=
\big(\sum_{l\in\mathbb{Z}}2^{lrs_{1}}\|\D_{l}u(t)\|_{L^{\rho}(L^{p})}^{r}\big)^{\frac{1}{r}}\,.$$
We then define the space $\widetilde{L}^{\rho}_{T}(B^{s_{1}}_{p,r})$ as the set of temperate distribution $u$ over
$(0,T)\times\R^{N}$ such that %$\lim_{q\rightarrow+\infty}S_{q}u=0$ in ${\cal S}^{'}((0,T)\times\R^{N})$
%and
$\|u\|_{\widetilde{L}^{\rho}_{T}(B^{s_{1}}_{p,r})}<+\infty$.
\end{definition}
We set $\widetilde{C}_{T}(\widetilde{B}^{s_{1}}_{p,r})=\widetilde{L}^{\infty}_{T}(\widetilde{B}^{s_{1}}_{p,r})\cap
{\cal C}([0,T],B^{s_{1}}_{p,r})$.
Let us emphasize that, according to Minkowski inequality, we have:
$$\|u\|_{\widetilde{L}^{\rho}_{T}(B^{s_{1}}_{p,r})}\leq\|u\|_{L^{\rho}_{T}(B^{s_{1}}_{p,r})}\;\;\mbox{if}\;\;r\geq\rho
,\;\;\;\|u\|_{\widetilde{L}^{\rho}_{T}(B^{s_{1}}_{p,r})}\geq\|u\|_{L^{\rho}_{T}(B^{s_{1}}_{p,r})}\;\;\mbox{if}\;\;r\leq\rho
.$$
\begin{remarka}
It is easy to generalize proposition \ref{produit1},
to $\widetilde{L}^{\rho}_{T}(B^{s_{1}}_{p,r})$ spaces. The indices $s_{1}$, $p$, $r$
behave just as in the stationary case whereas the time exponent $\rho$ behaves according to H\"older inequality.
\end{remarka}
In the sequel we will need of composition lemma in $\widetilde{L}^{\rho}_{T}(B^{s}_{p,r})$ spaces.
\begin{lemme}
\label{composition}
Let $s>0$, $(p,r)\in[1,+\infty]$ and $u\in \widetilde{L}^{\rho}_{T}(B^{s}_{p,r})\cap L^{\infty}_{T}(L^{\infty})$.
\begin{enumerate}
 \item Let $F\in W_{loc}^{[s]+2,\infty}(\R^{N})$ such that $F(0)=0$. Then $F(u)\in \widetilde{L}^{\rho}_{T}(B^{s}_{p,r})$. More precisely there exists a function $C$ depending only on $s$, $p$, $r$, $N$ and $F$ such that:
$$\|F(u)\|_{\widetilde{L}^{\rho}_{T}(B^{s}_{p,r})}\leq C(\|u\|_{L^{\infty}_{T}(L^{\infty})}\|u\|_{\widetilde{L}^{\rho}_{T}(B^{s}_{p,r})}.$$
\item Let $F\in W_{loc}^{[s]+3,\infty}(\R^{N})$ such that $F(0)=0$. Then $F(u)-F^{'}(0)u\in \widetilde{L}^{\rho}_{T}(B^{s}_{p,r})$. More precisely there exists a function $C$ depending only on $s$, $p$, $r$, $N$ and $F$ such that:
$$\|F(u)-F^{'}(0)u\|_{\widetilde{L}^{\rho}_{T}(B^{s}_{p,r})}\leq C(\|u\|_{L^{\infty}_{T}(L^{\infty})}\|u\|^{2}_{\widetilde{L}^{\rho}_{T}(B^{s}_{p,r})}.$$
\end{enumerate}
\end{lemme}
Here we recall a result of interpolation which explains the link
of the space $B^{s}_{p,1}$ with the space $B^{s}_{p,\infty}$, see
\cite{DFourier}.
\begin{proposition}
\label{interpolationlog}
There exists a constant $C$ such that for all $s\in\R$, $\e>0$ and
$1\leq p<+\infty$,
$$\|u\|_{\widetilde{L}_{T}^{\rho}(B^{s}_{p,1})}\leq C\frac{1+\e}{\e}\|u\|_{\widetilde{L}_{T}^{\rho}(B^{s}_{p,\infty})}
\biggl(1+\log\frac{\|u\|_{\widetilde{L}_{T}^{\rho}(B^{s+\e}_{p,\infty})}}
{\|u\|_{\widetilde{L}_{T}^{\rho}(B^{s}_{p,\infty})}}\biggl).$$ \label{5Yudov}
\end{proposition}
Now we give some result on the behavior of the Besov spaces via some pseudodifferential operator (see \cite{DFourier}).
\begin{definition}
Let $m\in\R$. A smooth function function $f:\R^{N}\rightarrow\R$ is said to be a ${\cal S}^{m}$ multiplier if for all muti-index $\alpha$, there exists a constant $C_{\alpha}$ such that:
$$\forall\xi\in\R^{N},\;\;|\p^{\alpha}f(\xi)|\leq C_{\alpha}(1+|\xi|)^{m-|\alpha|}.$$
\label{smoothf}
\end{definition}
\begin{proposition}
Let $m\in\R$ and $f$ be a ${\cal S}^{m}$ multiplier. Then for all $s\in\R$ and $1\leq p,r\leq+\infty$ the operator $f(D)$ is continuous from $B^{s}_{p,r}$ to $B^{s-m}_{p,r}$.
\label{singuliere}
\end{proposition}
%Actually, in \cite{BCD}, the proposition below is proved for non-homogeneous Besov spaces. The adaptation to homogeneous spaces is straightforward.
Let us now give some estimates for the heat equation:
\begin{proposition}
\label{5chaleur} Let $s\in\R$, $(p,r)\in[1,+\infty]^{2}$ and
$1\leq\rho_{2}\leq\rho_{1}\leq+\infty$. Assume that $u_{0}\in B^{s}_{p,r}$ and $f\in\widetilde{L}^{\rho_{2}}_{T}
(\widetilde{B}^{s-2+2/\rho_{2}}_{p,r})$.
Let u be a solution of:
$$
\begin{cases}
\begin{aligned}
&\p_{t}u-\mu\D u=f\\
&u_{t=0}=u_{0}\,.
\end{aligned}
\end{cases}
$$
Then there exists $C>0$ depending only on $N,\mu,\rho_{1}$ and
$\rho_{2}$ such that:
$$\|u\|_{\widetilde{L}^{\rho_{1}}_{T}(\widetilde{B}^{s+2/\rho_{1}}_{p,r})}\leq C\big(
 \|u_{0}\|_{B^{s}_{p,r}}+\mu^{\frac{1}{\rho_{2}}-1}\|f\|_{\widetilde{L}^{\rho_{2}}_{T}
 (\widetilde{B}^{s-2+2/\rho_{2}}_{p,r})}\big)\,.$$
 If in addition $r$ is finite then $u$ belongs to $C([0,T],B^{s}_{p,r})$.
\end{proposition}
\subsection{Hybrid Besov spaces}
The homogeneous Besov spaces fail to have nice inclusion properties: owing to the low frequencies, the embedding $B^{s}_{p,1}\hookrightarrow B^{t}_{p,1}$ does not hold for $s>t$. Still, the functions of $B^{s}_{p,1}$ are locally more regular than those of $B^{t}_{p,1}$: for any $\phi\in C^{\infty}_{0}$ and $u\in B^{s}_{p,1}$, the function $\phi u\in B^{t}_{p,1}$. This motivates the definition of Hybrid Besov spaces introduced by R. Danchin in \cite{DG} where the growth conditions satisfied by the dyadic blocks and the coefficient of integrability are not the same for low and high frequencies. Hybrid Besov spaces have been used in \cite{DG1} to prove global well-posedness for compressible gases in critical spaces. We generalize here a little bit the definition by distinguishing the coefficients of integrability.
\begin{definition}
 \label{def1.9}
Let $s,t\in\R$ and $(p,q)\in[1,+\infty]$. We set:
$$\|u\|_{\widetilde{B}^{s,t}_{p,q,1}}=\sum_{l\leq 0}2^{ls}\|\D_{l}u\|_{L^{p}}+\sum_{l\leq 0}2^{lt}\|\D_{l}u\|_{L^{q}}.$$
\end{definition}
\begin{notation}
 We will often use the following notation:
$$u_{BF}=\sum_{l\leq0}\D_{l}u\;\;\;\mbox{and}\;\;\;u_{HF}=\sum_{l>0}\D_{l}u.$$
\end{notation}
\begin{remarka}
 We have the following properties:
\begin{itemize}
 \item We have $\widetilde{B}^{s,s}_{p,p,1}=B^{s}_{p,1}$.
\item If $s_{1}\geq s_{3}$ and $s_{2}\geq s_{4}$ then $\widetilde{B}^{s_{3},s_{2}}_{p,q,1}\h \widetilde{B}^{s_{1},s_{4}}_{p,q,1}$.
\end{itemize}
\end{remarka}
We shall also make use of hybrid Besov-spaces. For them, one can prove results analoguous to proposition  \ref{produit1}, we refer to proposition \ref{hybrid} in the appendix.
\section{The mass conservation equation}
\label{section3}
We begin this section by recalling some estimates in Besov spaces for transport and heat equations. For more details, the reader is referred to \cite{BCD}.
\begin{proposition}
Let $1\leq p_{1}\leq p\leq+\infty$, $r\in[1,+\infty]$ and $s\in\R$ be such that:
$$-N\min(\frac{1}{p_{1}},\frac{1}{p^{'}})<s<1+\frac{N}{p_{1}}.$$
Suppose that $q_{0}\in B^{s}_{p,r}$, $F\in L^{1}(0,T, B^{s}_{p,r})$ and that $q\in
L^{\infty}_{T}(B^{s}_{p,r})\cap
C([0,T];{\cal S}^{'})$ solves the following transport equation:
$$
\begin{cases}
\begin{aligned}
&\p_{t}q+u\cdot\n q=F,\\
&q_{\ t=0}=q_{0}.
\end{aligned}
\end{cases}
$$
There exists a constant $C$ depending only on $N$, $p$, $p_{1}$, $r$ and $s$ such that , we have for a.e $t\in[0,T]$:
\begin{equation}
\|q\|_{\widetilde{L}^{\infty}_{t}(B^{s}_{p,r})}\leq e^{CU(t)}\big(\|q_{0}\|_{B^{s}_{p,r}}+\int^{t}_{0}e^{-CU(\tau)}
\|F(\tau)\|_{B^{s}_{p,r}}d\tau\big),
\label{20}
\end{equation}
with:
$U(t)=\int^{t}_{0}\|\n u(\tau)\|_{B^{\frac{N}{p_{1}}}_{p_{1},\infty}\cap L^{\infty}}d\tau$.
%If in addition $\D_{l}f\in C([0,T];L^{2})$ for all $l\in\mathbb{Z}$, and $r<+\infty$ then $f$ belongs to
%$\widetilde{C}_{T}(B^{s}_{2,r})$.
\label{transport1}
\end{proposition}
We want study now the following problem:
$$
\begin{cases}
\begin{aligned}
&\p_{t}q+u\cdot\n q+\alpha q=F,\\
&q_{\ t=0}=q_{0}.
\end{aligned}
\end{cases}
\leqno{({\cal H})}
$$
Above $a$ is the unknown function. We assume that $F\in L^{r}(0,T;B^{s}_{p,r})$, that $v$ is time dependent
vector-fields with coefficients in $L^{1}(0,T;B^{\frac{N}{p_{1}}+1}_{p_{1},1})$ and $\alpha>0$.
\\
Indeed we recall that we can rewrite the transport equation on the following form:
$$\p_{t}q+u\cdot\n q+(q+1)(P(1+q)-P(1))=-(1+q){\rm div}v_{1},$$
where we refer to the section \ref{section41} for the definition of $v_{1}$.
%, that $b$ is bounded by below by a positive constant $\underline{b}$ and %that
%$b$ belongs to $L^{\infty}(0,T;B^{\frac{N}{p}}_{p,1})$ with $p\in[1,+\infty]$.
\begin{proposition}
Let $1\leq p_{1}\leq p\leq+\infty$, $r\in[1,+\infty]$ and $s\in\R$ be such that:
$$-N\min(\frac{1}{p_{1}},\frac{1}{p^{'}})<s<1+\frac{N}{p_{1}}.$$
There exists a constant $C$ depending only on $N$, $p$, $p_{1}$, $r$ and $s$ such that for all $a\in L^{\infty}([0,T],B^{\sigma}_{p,r})$ of $({\cal H})$ with initial data $a_{0}$ in $B^{s}_{p,r}$ and $g\in L^{1}([0,T], B^{s}_{p,r})$, we have for a.e $t\in[0,T]$:
\begin{equation}
\|q\|_{\widetilde{L}^{\infty}_{t}(B^{s}_{p,r})}+\|q\|_{\widetilde{L}^{1}_{t}(B^{s}_{p,r})}\leq e^{CU(t)}\big(\|q_{0}\|_{B^{s}_{p,r}}+\int^{t}_{0}e^{-CU(\tau)}
\|F(\tau)\|_{B^{s}_{p,r}}d\tau\big),
\label{20}
\end{equation}
with:
$U(t)=\int^{t}_{0}\|\n u(\tau)\|_{B^{\frac{N}{p_{1}}}_{p_{1},\infty}\cap L^{\infty}}d\tau$.
\label{transport2}
\end{proposition}
{\bf Proof:}\;\;%{\bf Proof:}
Applying $\D_{l}$ to $({\cal H})$ yields:
$$\p_{t}\D_{l}q+u\cdot\n\D_{l}q+\alpha\D_{l}q=R_{l}%+\widetilde{R}_{l}
+\D_{l}F,$$
with $R_{l}=[u\cdot\n,\D_{l}]q$. %and $\widetilde{R}_{l}=[b,\D_{l}]a$.
Multiplying by $\D_{l}a|\D_{l}a|^{p-2}$ then performing a time integration, we easily get:
$$
\begin{aligned}
&\|\D_{l}q(t)\|_{L^{p}}+\alpha\int^{t}_{0}\|\D_{l}q(s)\|_{L^{p}}ds\leq\|\D_{l}q_{0}\|_{L^{p}}
+\int^{t}_{0}\big(\|R_{l}\|_{L^{p}}+\frac{1}{p}\|{\rm div}u\|_{L^{\infty}}\|\D_{l}q\|_{L^{p}}\\
&\hspace{11,5cm}+\|\D_{l}F\|_{L^{p}}\big)d\tau.
\end{aligned}
$$
Next the term $\|R_{l}\|_{L^{p}}$% and $\|\widetilde{R}_{l}\|_{L^{p}}$
may be bounded according to lemma \ref{alemme2} in appendix. We get then:
$$
\begin{aligned}
&\|q\|_{\widetilde{L}^{\infty}_{t}(B^{s}_{p,r})}+\alpha\|q\|_{\widetilde{L}^{1}_{t}(B^{s}_{p,r})}ds\leq\|\D_{l}q_{0}\|_{B^{s}_{p,r}}
+\int^{t}_{0}\big(\|F(\tau)\|_{B^{s}_{p,r}}+C U^{'}(\tau)\|q\|_{\widetilde{L}^{\infty}_{t}(B^{s}_{p,r})}\big)d\tau.
\end{aligned}
$$
We end up with Gronwall lemma by letting $X(t)=\|q\|_{\widetilde{L}^{\infty}_{t}(B^{s}_{p,r})}+\alpha\|q\|_{\widetilde{L}^{1}_{t}(B^{s}_{p,r})}$.
%\begin{equation}
%\forall t\in [0,T],\;\forall l\in\mathbb{Z},\;\;2^{ls}\|\D_{l}a(t)\|_{L^{p}}\leq 2^{ls}\|\D_{l}a_{0}\|_{L^{p}}+
%C\int^{t}_{0}c_{l}(1+\|a\|_{B^{s}_{p,r}\cap L^{\infty}})V^{'}d\tau,
%\label{25}
%\end{equation}
%hence, summing up on $\mathbb{Z}$ in $l^{r}$,
%$$\forall t\in [0,T],\;\forall %l\in\mathbb{Z},\;\;\|a(t)\|_{B^{s}_{p,r}}\leq\|a_{0}\|_{B^{s}_{p,r}}+\int^{t}_{0}CV^{'}\|a(\tau)\|_{B^{s}_{p,r}}d\tau
%+\int^{t}_{0}C(1+\|a\|_{L^{\infty}_{T}})V^{'}d\tau.
%$$
%Next we have:
%$$\|a\|_{L^{\infty}_{t}}\leq\int^{t}_{0}(1+\|a(\tau)\|_{L^{\infty}})V^{'}(\tau)d\tau.$$
%By summing the two previous inequalities, applying Gronwall lemma and proposition \ref{resteimp1} yields inequality (\ref{22}).
\section{The proof of theorem \ref{theo1}}
\label{section4}
\subsection{Strategy of the proof}
\label{section41}
To improve the results of Danchin in \cite{DG}, Charve and Danchin in \cite{CD} and Chen et al in \cite{CMZ}, it is crucial to kill the coupling between the velocity and the pressure which exists in these works . In this goal, we need to integrate the pressure term in the study of the linearized equation of the momentum equation as in \cite{H2}. For making, we will try to express the gradient of the pressure as a Laplacian term, so we set for $\bar{\rho}>0$ a constant state:
$${\rm div}v=P(\rho)-P(\bar{\rho}).$$
Let ${\cal E}$ the fundamental solution of the Laplace operator.
$$$$
We will set in the sequel: $v=\n{\cal E}*\big(P(\rho)-P(\bar{\rho})\big)=\n\big({\cal E}*[P(\rho)-P(\bar{\rho})]\big)$ ( $*$ here means the operator of convolution). We verify next that:
$$
\begin{aligned}
\n{\rm div}v=\n\D \big({\cal E}*[P(\rho)-P(\bar{\rho})]\big)=\D\n\big({\cal E}*[P(\rho)-P(\bar{\rho})]\big)=\D v=\n P(\rho).
\end{aligned}
$$
%So clearly, we have $\n{\rm div}v=\n P(\rho)$.
%Next by using the Helmholtz decomposition, we can writte $v$ as a sum of a free divergence field $h$ and of a gradient vector $\n F$ as follows
%$v=\n F+h$
%where ${\rm div}h=0$.
%Formally we will set in the sequel $v$ and $F$ the unique solution of the following system:
%$$
%\begin{cases}
%\begin{aligned}
%&{\rm div}v=P(\rho)-P(\bar{\rho}),\\
%&v=\n F.
%\end{aligned}
%\end{cases}
%$$
%We have then $\D F= P(\rho)-P(\bar{\rho})$ and $\D\n F=\n P(\rho)=\D v$. In the sequel we will assume always that $h=0$. By this way we have shown that:
%$$\D v=\n{\rm div}v=\n P(\rho).$$
By this way we can now rewrite the momentum equation of (\ref{0.6}). We obtain the following equation where we have set $\nu=2\mu+\lambda$:
$$\p_{t}u+u\cdot \n u-\frac{\mu}{\rho}\D\big(u-\frac{1}{\nu}v\big)-\frac{\lambda+\mu}{\rho}\n{\rm div}\big(u-\frac{1}{\nu}v\big)=f.$$
We want now calculate $\p_{t}v$, by the transport equation we get:
$$\p_{t}v=\n{\cal E}*\p_{t}P(\rho)=-\n {\cal E}*\big(P^{'}(\rho){\rm div}(\rho u)\big).$$
We have finally:
$$\D(\p_{t}F)=-P^{'}(\rho){\rm div}(\rho u).$$
\begin{notation}
To simplify the notation, we will note in the sequel
$$\n {\cal E}*\big(P^{'}(\rho){\rm div}(\rho u)\big)=\n(\D)^{-1}\big(P^{'}(\rho){\rm div}(\rho u)\big).$$
\end{notation}
Finally we can now rewrite the system (\ref{0.6}) as follows:
\begin{equation}
\begin{cases}
\begin{aligned}
&\p_{t}q+(v_{1}+\frac{1}{\nu}v)\cdot\n q+\frac{1}{\nu}(1+q)(P(\rho)-P(1))=-(1+q){\rm div}v_{1},\\
&\p_{t}v_{1}-\frac{1}{1+q}{\cal A}v_{1}=f-u\cdot\n u+\frac{1}{\nu}\n(\D)^{-1}\big(P^{'}(\rho){\rm div}(\rho u)\big),\\
&q_{/ t=0}=a_{0},\;(v_{1})_{/ t=0}=(v_{1})_{0}.
\end{aligned}
\end{cases}
\label{0.7}
\end{equation}
where $v_{1}=u-\frac{1}{\nu}v$ is called the effective velocity. In the sequel we will study this system by exctracting some uniform bounds in Besov spaces on
$(q,v_{1})$. %as the in the following works \cite{AP}, \cite{H}.
The advantage of the system (\ref{0.7}) is that we have \textit{kill} the coupling between $v_{1}$ and a term of pressure. Indeed in the works \cite{CD} and \cite{CMZ}, the pressure was included in the study of the linear system, it means mandatory a coupling between the density and the velocity. In particular it was impossible to distinguish the index of integration for the Besov spaces.
\\
However we can remark that this change of variable $v_{1}$ is interesting only in the case of low frequencies, indeed heuristically in low frequencies $\n P(\rho)$ is small in Fourier variable so it is not a matter.\\
It is natural in this case to study the variable  $u$ in low frequencies. Moreover as explained in the introduction, the system (\ref{0.1}) has a hyperbolic behavior, which means that we can work only with spaces builtet on $L^{2}$ (indeed classicaly the hyperbolic system are ill-posed in space constructed on general $L^{p}$). In the following ssection, we will explain how to treat the case in low frequencies and how we will use the fact that $q$ behaves in low frequencies as an heat equation.
%$$
%\begin{aligned}
%&\p_{t}v_{1}+v_{1}\cdot\n v_{1}-\frac{\mu}{\rho}\D v_{1}-\frac{\lambda+\mu}{\rho}\n{\rm div}v_{1}=f-\frac{1}{\nu^{2}}v\cdot\n v\\
%&\hspace{4cm}-\frac{1}{\nu}(v\cdot\n v_{1}+v_{1}\cdot\n v)+\frac{1}{\nu}\n(\D)^{-1}\big(P^{'}(\rho){\rm div}(\rho u)\big).
%\end{aligned}
%$$
\subsection{A linear model with convection}
In this section, we will explain how we treat the low frequency regime by following the approch of Charve and Danchin in \cite{CD}.
In low frequencies, the first order terms predominate and the viscous term $\D u$ may be neglected so that (\ref{0.1}) has to be treated by means of hyperbolic energy methods. It means that we can only work in spaces constructed on $L^{2}$. Moreover in the case of low frequencies the effective velocity is not a adapted variable in the sense that it is less regular than $u$ as $(\D)^{-1}\n P(\rho)$ is not very regular. It is better in this case to work with $u$.%we give the sketch of the proof of theorem
%\ref{5theo1} on the global existence result with small initial
%data.
The first idea would be to study the linear system associated to
(\ref{0.1}), it means:%. We concentrate on the first two equations because the third
%equation is just a heat equation with a non linear term. The system
%we want to study reads:
$$
\begin{cases}
\begin{aligned}
&\p_{t}q+{\rm div}u=F^{'},\\
&\p_{t}u-\mu \D u-\lambda\n {\rm div}u+\n q=G^{'}.
\end{aligned}
\end{cases}
\leqno{(PH)}
$$
This system has been studied by D. Hoff and K. Zumbrum in \cite{5HZ}. There, they investigate the decay
estimates, and exhibit the parabolic smoothing effect on $u$ and on
the low frequencies of
$q$, and a damping effect on the high frequencies of $q$.\\
The problem is that if we focus on this linear system, it appears
impossible to control the term of convection $u\cdot\n q$ which is
one derivative less regular than $q$. However in low frequencies the Green matrix of the linearized systems behaves as the heat kernel
(see \cite{CMZ}), the terms $v\cdot\n q$ and $v\cdot\n u$ can be handled as the perturbation terms. We study then the following system:
%.Hence we shall include the
%convection term in the linear system. We thus have to study:
$$
\begin{cases}
\begin{aligned}
&\p_{t}q+{\rm div}u=-+v\cdot\n q+F,\\
&\p_{t}d-\mu \D u-\lambda\n {\rm div}u+\n q=-v\cdot\n u+G,
\end{aligned}
\end{cases}
\leqno{(LH)^{'}}
$$
%where $v$ is a function  and we will precise its regularity in the
%next proposition. In fact% R. Danchin in \cite{DG}, and
%F. Charve and R. Danchin in \cite{CD} have obtained estimates in Besov-type space for the following paralinerization of system (\ref{0.1}):
%$$
%\begin{cases}
%\begin{aligned}
%&\p_{t}q+{\rm div}(T_{v}q)+{\rm div}u=F,\\
%&\p_{t}u+T_{v}\cdot\n u--\mu \D u-\lambda\n {\rm div}u+\n q=G,
%\end{aligned}
%\end{cases}
%\leqno{(PL)}
%$$
%with ${\rm div}(T_{v}q)=\p_{i}(T_{v}q)$ and $T_{v}\cdot\n u=T_{v_{i}}\p_{i}u$. We recall that $T_{v}u==\sum_{q}S_{q-1}u\D_{q}V$.
%Here we want point out that the fact to include ${\rm div}u$ in the linear system is crucial, indeed it we we see this term as a term of rest we lost one derivative. Let us first recall a priori estimate for $(PL)$ in Besov spaces moded on $L^{2}$, we refer to proposition 6 in \cite{CD}.
%System $(SW2)^{'}$ has been
%studied  by R. Danchin in \cite{DG}, he obtains the following proposition.
We obtain then the following proposition:
\begin{proposition}
\label{5linear1}
Let $(q,u)$ a solution of $(LH)^{'}$, let $s\in \R$. % with $1-\N<s\leq 1+\N$
 %and $V(t)=\int^{t}_{0}\|\n v(\tau)\|_{L^{\infty}}d\tau$. We assume that $(F,G)_{HF}=0$.
The following estimate holds:
$$
\begin{aligned}
&\|(q,u)_{BF}\|_{\widetilde{L}^{\infty}(\widetilde{B}^{s}_{2,1})}+\|(q,u)_{BF}\|_{\widetilde{L}^{1}(\widetilde{B}^{s+2}_{2,1})}\leq
%+\|d(t)\|_{B^{s-1}_{2,1}}+\bar{\nu}\int^{t}_{0}(\|q(\tau)\|_{\widetilde{B}^{s+1,s}_{2,1}}+
%\|d(\tau)\|_{B^{s+1}_{2,1}})d\tau
%\leq\\% e^{CV(t)}\\
\|(q_{0},u_{0})_{BF}\|_{\widetilde{B}^{s}_{2,1}}+\|(F,G)_{BF}\|_{\widetilde{L}^{1}(\widetilde{B}^{s}_{2,1})}\\
&\hspace{9cm}+\|(v\cdot\n q,v\cdot\n u){BF}\|_{\widetilde{L}^{1}(\widetilde{B}^{s}_{2,1})}.
\end{aligned}
$$
\end{proposition}
{\bf Proof:}\\
\\
In this case for $j\leq 0$, in terms of Green matrix (see \cite{CMZ}), the solution of $(LH)^{'}$ can be expressed as:
$$\left(\begin{array}{c}
\D_{j}q(t)\\
\D_{j} u(t)\\
\end{array}
\right) =W(t)\left(\begin{array}{c}
\D_{j}q_{0}\\
\D_{j}u_{0}\\
\end{array}
\right)+\int_{0}^{t}W(t-s)\left(\begin{array}{c}
\D_{j}F(s)-\D_{j}(v\cdot\n q)\\
\D_{j}G(s)-\D_{j}(v\cdot\n u)\\
\end{array}
\right)\ ds\;.$$
with $W$ the Green matrix.
From proposition 4.4 in \cite{CMZ} and Young's inequality we obtain the result.{\hfill $\Box$}
\subsection{Proof of the existence}
\subsubsection*{Construction of approximate solutions}
We use a standard scheme:
\begin{enumerate}
\item We smooth out the data and get a sequence of global smooth solutions $(q^{n},u^{n})_{n\in\mathbb{N}}$ to (\ref{0.1})
on $\R$  by using the results of \cite{CD} and \cite{CMZ}.% which may depend on $n$. We set $v_{1}^{n}=u^{n}-\frac{1}{\nu}v^{n}$ where ${\rm div}v^{n}=P(\rho^{n})-P(\bar{\rho})$.
\item %We exhibit a positive lower bound $T$ for $T^{n}$, and
We prove uniform estimates on $(q^{n},v_{1}^{n})$ in high frequencies and on $(q^{n},u^{n})$ in low frequencies. %the space
%$$E=\widetilde{C}(B^{\NN}_{p,1}\cap B^{\NN-1}_{p,1})\times\big(\widetilde{C}_{T}(B^{\frac{N}{p_{1}}-1}_{p_{1},1}+B^{\NN-1}_{p,1})\cap\widetilde{L}^{1}_{T}(
%B^{\frac{N}{p_{1}}+1}_{p_{1},1}+B^{\NN+1}_{p,1})\big).$$
%and we will deduce uniform estimates on $(q^{n},u^{n})$.
%More precisely to get this bounds we will need to study the behavior of $(a^{n},v_{1}^{n})$.
%for the smooth solution $(a^{n},v_{1}^{n})$.
\item We use compactness to prove that the sequence $(q^{n},u^{n})$ converges, up to extraction, to a solution of (\ref{0.1}).
\end{enumerate}
%Throughout the proof, we denote $\underline{\nu}=\underline{b}\min(\mu,\lambda+2\mu)$ and $\bar{\nu}=\mu+|\mu+\lambda|$, and we assume (with no loss of generality) that $f$ belongs to $\widetilde{L}^{1}_{T}(B^{\frac{N}{p}-1}_{p_{1},1})$.
\subsubsection*{First step}
We smooth out the data as follows:
$$q_{0}^{n}=S_{n}q_{0},\;\;u_{0}^{n}=S_{n}u_{0}\;\;\;\mbox{and}\;\;\;f^{n}=S_{n}f.$$
Note that we have:
$$\forall l\in\mathbb{Z},\;\;\|\D_{l}q^{n}_{0}\|_{L^{p}}\leq\|\D_{l}q_{0}\|_{L^{p}}\;\;\;\mbox{and}\;\;\;\|q^{n}_{0}\|
_{B^{\frac{N}{p}}_{p,1}\cap B^{\NN-1}_{p,1}}\leq \|q_{0}\|_{B^{\frac{N}{p}}_{p,1}\cap B^{\NN-1}_{p,1}},$$
and similar properties for $u_{0}^{n}$ and $f^{n}$, a fact which will be used repeatedly during the next
steps. Now, according \cite{DG}, one can solve (\ref{0.1}) with the smooth data $(q_{0}^{n},u_{0}^{n},f^{n})$.
We get a solution $(q^{n},u^{n})$ %on a non trivial time interval $[0,T_{n}]$
such that:
\begin{equation}
\begin{aligned}
&q^{n}\in\widetilde{C}(\R,B^{N}_{2,1}\cap B^{\N-1}_{2,1})\;\;\mbox{and}\;\;u^{n}\in\widetilde{C}(\R,B^{\N-1}_{2,1})\cap
\widetilde{L}^{1}
(\R,B^{\N+1}_{2,1}).
\end{aligned}
\label{a26}
\end{equation}
\subsubsection*{Uniform bounds}
%Let $T_{n}$ be the lifespan of $(a_{n},u_{n})$, that is the supremum of all $T>0$ such that (\ref{0.1}) with initial data
%$(a_{0}^{n},u_{0}^{n})$ has a solution which satisfies (\ref{a26}). Let $T$ be in $(0,T_{n})$.
%$u^{n}=u_{L}+\bar{u}^{n}$ with:
%$$
%\begin{cases}
%\begin{aligned}
%&\p_{t}u_{L}-\mu\D u_{L}-(\mu+\lambda)\n{\rm div}u_{L}=0,\\
%&(u_{L})_{t=0}=u_{0},
%\end{aligned}
%\end{cases}
%$$
We set now
$$v_{n}=\n\big({\cal E}*[P(\rho^{n})-P(1)]\big)\;\;\;\mbox{with}\;\;\;{\rm div}v^{n}=P(\rho^{n})-P(1)\;\;\;\mbox{and}\;\;\;v_{1}^{n}=u^{n}-\frac{1}{\nu}v^{n},$$
with ${\cal E}$ the fundamental solution of the Laplace operator and $\nu=\lambda+\mu$.
In the sequel we will note $g(q^{n})=P(\rho^{n})-P(1)$ where $g$ is a regular function.
In this part, we aim at getting uniform estimates on $(q^{n}_{HF},(v_{1})_{HF}^{n})$ in high frequencies and on $(q^{n}_{BF},u^{n}_{BF})$ in low frequencies in the following space $E^{'}$ and $F^{'}$:
$$
\begin{aligned}
 &E^{'}=\big(\widetilde{L}^{\infty}(B^{\NN}_{p,1})\cap\widetilde{L}^{1}(B^{\NN}_{p,1})\big)\times\big(\widetilde{L}^{\infty}
(B^{\frac{N}{p_{1}}-1}_{p_{1},1}+B^{\NN}_{p,1})+\widetilde{L}^{1}(B^{\frac{N}{p_{1}}+1}_{p_{1},1}+B^{\NN+2}_{p,1})\big).\\
&F^{'}=\big(\widetilde{L}^{\infty}(B^{\N-1}_{2,1})\cap\widetilde{L}^{1}(B^{\N+1}_{2,1})\big)\times\big(\widetilde{L}^{\infty}
(B^{\frac{N}{2}-1}_{2,1})+\widetilde{L}^{1}(B^{\frac{N}{2}+1}_{2,1})\big).\\
\end{aligned}
$$
More precisely we will obtain uniform estimates on $(q^{n},u^{n})$ in $E$ and on $(q^{n},v_{1}^{n})$ in $F$ whith:
$$
\begin{aligned}
E=&\big(\widetilde{L}^{\infty}(\widetilde{B}^{\N-1,\NN}_{2,p,1})\cap \widetilde{L}^{1}(\widetilde{B}^{\N+1,\NN}_{2,p,1})\big)\times
\big(\widetilde{L}^{\infty}(\widetilde{B}^{\N-1,\NN}_{2,p,1}+\widetilde{B}^{\N-1,\frac{N}{p_{1}}-1}_{2,p_{1},1})\\
&\hspace{8cm}\cap \widetilde{L}^{1}(\widetilde{B}^{\N+1,\NN+1}_{2,p,1})\big).
\end{aligned}
$$
$$
\begin{aligned}
F=&\widetilde{L}^{\infty}(\widetilde{B}^{\N-1,\NN}_{2,p,1})\cap \widetilde{L}^{1}(\widetilde{B}^{\N+1,\NN}_{2,p,1})\times
\big(\widetilde{L}^{\infty}(\widetilde{B}^{\N,\NN}_{2,p,1}+\widetilde{B}^{\N,\frac{N}{p_{1}}-1}_{2,p_{1},1})\\
&\hspace{6cm}\cap \widetilde{L}^{1}(\widetilde{B}^{\N+2,\NN+2}_{2,p,1}+\widetilde{B}^{\N+2,\frac{N}{p_{1}}+1}_{2,p_{1},1})\big).
\end{aligned}
$$
We will work finally in the space $H$ with:
$$(q,u)\in H\Leftrightarrow (q,u)_{BF}\in E^{'}\;\;\mbox{and}\;\;(q,v_{1})_{HF}\in F^{'}.$$
We have then: $\|(q,u)\|_{H}=\|(q,u)_{BF}\|_{E^{'}}+\|(q,v_{1})_{HF}\|_{F^{'}}$.
% for all $T<T_{n}$. This will show
%that $T_{n}=+\infty$ by continuity arguments.
We can now check that $(q^{n},v^{n}_{1})$ satisfy the following system:
\begin{equation}
\begin{cases}
\begin{aligned}
&\p_{t}q^{n}+u^{n}\cdot\n q^{n}+\frac{P^{'}(1)}{\nu}q^{n}=F^{n}_{1},\\
&\p_{t}v_{1}^{n}-{\cal A}v_{1}^{n}=F^{n}_{2}+f,\\
&q^{n}_{0}=q_{0},\;(v_{1}^{n})_{/ t=0}=u_{0}^{n}-\frac{1}{\nu}v_{0}^{n}.
%-\frac{1}{\nu}\widetilde{v}^{n}_{0}.
\end{aligned}
\end{cases}
\label{systemessen}
\end{equation}
which is a transport equation and a heat equation.%has been studied in proposition \ref{5linear1} and with:
$$
\begin{aligned}
&F_{1}^{n}=-(1+q^{n}){\rm div}v_{1}^{n}-\frac{1}{\nu}(P(1+q^{n})-P(1)-P^{'}(1)q^{n})-\frac{1}{\nu}q^{n}(P(1+q^{n})-P(1)),\\
&G_{1}^{n}=(\frac{1}{1+q^{n}}-1)
{\cal A}v_{1}^{n}-u^{n}\cdot\n u^{n}+\frac{1}{\nu}\n(\D)^{-1}(P^{'}(\rho^{n}){\rm div}(\rho^{n}u^{n})).
\end{aligned}
$$
%In the sequel we have to decompose $\widetilde{v}_{1}^{n}$ due
%Define $m\in\mathbb{Z}$ by:
%\begin{equation}
%m=\inf\{ p\in\mathbb{Z}/\;2\bar{\nu}\sum_{l\geq p}2^{l\frac{N}{p}}\|\D_{l}a_{0}\|_{L^{p}}\leq c\bar{\nu}\}
%\label{def}
%\end{equation}
%where $c$ is small enough positive constant (depending only $N$) to be fixed hereafter. In the sequel we will need of a control on $a-S_{m}a$ small to apply proposition \ref{linearise}, so here $m$ is enough big (we explain how in the sequel).
Moreover $(q^{n},u^{n})_{n\in\mathbb{N}}$ is the solution of the following system:
\begin{equation}
\begin{cases}
\begin{aligned}
&\p_{t}q^{n}+u^{n}\cdot\n q^{n}+{\rm div}u^{n}=F^{n}\\
&\p_{t}u^{n}+u^{n}\cdot\n u^{n}-{\cal A}u^{n}+P^{'}(1)\n q^{n}=G^{n}+f^{n}\\
&(q^{n},u^{n})_{/t=0}=(q^{n}_{0},u^{n}_{0}),
\end{aligned}
\end{cases}
\label{5grandeequation}
\end{equation}
which has been studied for low frequencies in proposition \ref{5linear1} with:
$$
\begin{aligned}
&F^{n}=-q^{n}{\rm div}u^{n},\\
&G^{n}=-\frac{q^{n}}{1+q^{n}}{\cal A}u^{n}+(P^{'}(1)-P^{'}(1+q^{n}))\n q^{n}.
\end{aligned}
$$
Let us set:
$$
\begin{aligned}
%&\beta=\|u_{0}\|_{\widetilde{B}^{\N-1,\frac{N}{p_{1}}-1}_{2,p_{1},1}}+
%\|q_{0}\|_{\widetilde{B}^{\N-1,\frac{N}{p}}_{2,p,1}}+
%\|f\|_{L^{1}_{T}(B^{\frac{N}{p_{1}}-1}_{p_{1},1})},\\
&E(q,u)=\|q\|_{\widetilde{L}^{\infty}(\widetilde{B}^{\N-1,\NN}_{2,p,1})}+\|u\|_{\widetilde{L}^{\infty}(\widetilde{B}^{\N-1, \frac{N}{p_{1}}-1}_{2,p_{1},1}+\widetilde{B}^{\N-1, \frac{N}{p}}_{2,p,1})}+\|q\|_{\widetilde{L}^{1}(\widetilde{B}^{\N+1,\NN}_{2,p,1})}\\
&\hspace{8cm}+
\|u\|_{\widetilde{L}^{1}(\widetilde{B}^{\N+1,\frac{N}{p}+1}_{2,p,1})},\\
&E_{1}(q,u)=\|q\|_{\widetilde{L}^{\infty}(B^{\N-1}_{2,1})}+\|u\|_{\widetilde{L}^{\infty}(B^{\N-1}_{2,1})}+
\|q\|_{\widetilde{L}^{1}(B^{\N+1}_{2,1})}+
\|u\|_{\widetilde{L}^{1}(B^{\N+1}_{2,1})}.\\
&E_{2}(q,u)=\|q\|_{\widetilde{L}^{\infty}(B^{\NN}_{p,1})}+\|u\|_{\widetilde{L}^{\infty}(B^{ \frac{N}{p_{1}}-1}_{p_{1},1}+B^{ \frac{N}{p}}_{p,1})}+\|q\|_{\widetilde{L}^{1}(B^{\NN}_{p,1})}\\
&\hspace{8cm}+
\|u\|_{\widetilde{L}^{1}(B^{\frac{N}{p_{1}}+1}_{p_{1},1}+B^{\frac{N}{p}+2}_{p,1})}.
%\|q\|_{\widetilde{L}_{t}^{\infty}(\widetilde{B}^{\NN-1,\NN}_{p,1})}+\|q\|_{\widetilde{L}_{t}^{1}(\widetilde{B}^{\NN-1,\NN}_{p,1})}
%+\|v_{1}\|_{\widetilde{L}_{t}^{\infty}(B^{\NN-1}_{p,1}+B^{\frac{N}{p_{1}}-1}_{p_{1},1})}\\
%&\hspace{10cm}+\|v_{1}\|_{\widetilde{L}^{1}_{t}(
%B^{\NN+1}_{p,1}+B^{\frac{N}{p_{1}}+1}_{p_{1},1})}.
\end{aligned}
$$
One can now apply the propositions \ref{5chaleur}  at our system to obtain
uniform bounds, so %by setting $V_{n}(t)=\|u^{n}\|_{L^{1}_{T}(B^{\N+1})}$
we have in high frequencies to control $(v_{1}^{n},q^{n})$ and in low frequencies $(q^{n},u^{n})$:
$$
\begin{aligned}
&E_{2}((q^{n},v_{1}^{n})_{HF})\leq C\big(\|(q_{0})_{HF}\|_{B^{\NN-1}_{p,1}+B^{\NN}_{p,1}}+\|(u_{0})_{HF}\|_{B^{\frac{N}{p_{1}}-1}_{p_{1},1}}
\\
&\hspace{4cm}+\|(F_{1}^{n})_{HF}\|_{\widetilde{L}^{1}(B^{\NN}_{p,1})}
+\|G_{1}^{n}\|_{\widetilde{L}^{1}(B^{\frac{N}{p_{1}}-1}_{p_{1},1}+B^{\frac{N}{p}}_{p,1})}\big),
%\|v^{1}_{0}\|_{B^{\frac{N}{p_{1}}-1}_{p_{1},1}+B^{\NN-1}_{p,1}}+
%\|f\|_{\widetilde{L}^{1}(B^{\frac{N}{p_{1}}-1}_{p_{1},1})}+\|F^{n}_{2}\|_{\widetilde{L}^{1}(B^{\NN-1}_{p,1}+B^{\frac{N}{p_{1}}-1}_{p_{1}%,1})}\\
%&\hspace{4,5cm}+
%e^{CU^{n}(+\infty)}\big(\,\|q_{0}\|_{\widetilde{B}^{\NN-1,\NN}_{p,1}}
%+\int_{\R}e^{-CU_{n}(\tau)}\|F_{1}^{n}(\tau)\|_{\widetilde{B}^{\NN-1,\NN}_{p,1}}
%d\tau
%\big),
\end{aligned}
$$
and
$$
\begin{aligned}
&E_{1}((q^{n},u^{n})_{BF})\leq C\big(\|(q_{0})_{BF}\|_{B^{\N-1}_{2,1}}+\|(u_{0})_{BF}\|_{B^{\N-1}_{2,1}}
\\
&\hspace{4cm}+\|(F^{n})_{BF}\|_{\widetilde{L}^{1}(B^{\N-1}_{2,1})}
+\|G^{n}\|_{\widetilde{L}^{1}(B^{\N-1}_{2,1})}\big),
%\|v^{1}_{0}\|_{B^{\frac{N}{p_{1}}-1}_{p_{1},1}+B^{\NN-1}_{p,1}}+
%\|f\|_{\widetilde{L}^{1}(B^{\frac{N}{p_{1}}-1}_{p_{1},1})}+\|F^{n}_{2}\|_{\widetilde{L}^{1}(B^{\NN-1}_{p,1}+B^{\frac{N}{p_{1}}-1}_{p_{1}%,1})}\\
%&\hspace{4,5cm}+
%e^{CU^{n}(+\infty)}\big(\,\|q_{0}\|_{\widetilde{B}^{\NN-1,\NN}_{p,1}}
%+\int_{\R}e^{-CU_{n}(\tau)}\|F_{1}^{n}(\tau)\|_{\widetilde{B}^{\NN-1,\NN}_{p,1}}
%d\tau
%\big),
\end{aligned}
$$
%with:
%$$
%\begin{aligned}
%&F_{1}^{n}(\tau)=-(1+q^{n}){\rm div}v_{1}^{n}-\frac{1}{\nu}\big(P^{'}(1)q^{n}-(1+q^{n})g(q^{n})\big),\\
%&F_{2}^{n}(\tau)=(\frac{1}{1+q^{n}}-1)
%{\cal A}v_{1}^{n}-u^{n}\cdot\n u^{n}+\frac{1}{\nu}\n(\D)^{-1}(P^{'}(\rho^{n}){\rm div}(\rho^{n}u^{n})),\\
%&U_{n}(T)=\int^{T}_{0}\|\n u^{n}(\tau)\|_{B^{\NN}_{p,1}}d\tau,
%\end{aligned}
%$$
%where $g(q^{n})=P(\rho^{n}-P(1)$.
Therefore, it is only a matter of proving appropriate estimates for $F_{1}^{n}$, $G_{1}^{n}$,
$F^{n}$ and $G^{n}$  by using properties
of continuity on the paraproduct and proposition \ref{transport}, \ref{5chaleur} and \ref{5linear1}.\\
We begin by estimating
$\|(F_{1}^{n})_{HF}\|_{\widetilde{L}^{1}(B^{\NN}_{p,1})}$ and $
\|(G_{1}^{n})_{HF}\|_{\widetilde{L}^{1}(B^{\frac{N}{p_{1}}-1}_{p_{1},1}+B^{\frac{N}{p}}_{p,1})}$
, we have to use proposition \ref{produit1} and proposition \ref{hybrid} and the fact that by interpolation ${\rm div}v_{1}^{n}$
is in $\widetilde{L}^{1}(B^{\N+1,\NN}_{2,p,1})$ because $\widetilde{L}^{1}(\widetilde{B}^{\N+1,\NN+1}_{2,p,1}+\widetilde{B}^{\N+1,\frac{N}{p_{1}}}_{2,p_{1},1})\h \widetilde{L}^{1}(B^{\N+1,\NN}_{2,p,1})$ as $p_{1}\leq p$:
$$
\begin{aligned}
 \|\big((1+q^{n}){\rm div}v_{1}^{n}\big)_{HF}\|_{\widetilde{L}^{1}(B^{\NN}_{p,1})}&\leq
\|{\rm div}v_{1}^{n}\|_{\widetilde{L}^{1}(B^{\N+1,\NN}_{2,p,1})}+\|q\|_{L^{\infty}(L^{\infty})}\|{\rm div}v_{1}^{n}\|_{\widetilde{L}^{1}(B^{\N+1,\NN}_{2,p,1})}\\
&\hspace{3cm}+\|{\rm div}v_{1}^{n}\|_{L^{1}(L^{\infty})}\|q^{n}\|_{\widetilde{L}^{\infty}(B^{\N-1,\NN}_{2,p,1})}.
\end{aligned}
$$
$$
\begin{aligned}
\|\big[(P(1+q^{n})-P(1)-P^{'}(1)q^{n})\big]_{HF}\|_{\widetilde{L}^{1}(B^{\NN}_{p,1})}\leq C\|q^{n}\|^{2}_{\widetilde{L}^{2}(B^{\N,\NN}_{2,p,1})},
\end{aligned}
$$
$$\|\big[q^{n}(P(1+q^{n})-P(1))\big]_{HF}\|_{\widetilde{L}^{1}(B^{\NN}_{p,1})}\leq C\|q^{n}\|^{2}_{\widetilde{L}^{2}(B^{\N,\NN}_{2,p,1})},$$
Next we have to treat the term $[\frac{q^{n}}{1+q^{n}}{\cal A}v_{1}^{n}]_{HF}$ in $\widetilde{L}^{1}(B^{\frac{N}{p_{1}}-1}_{p_{1},1}+B^{\NN}_{p,1})$, where we can split ${\cal A}v_{1}^{n}$ on the form:
$$v_{1}^{n}=h^{n}+g^{n},$$
with: $h^{n}\in \widetilde{L}^{\infty}(B^{\N-1,\frac{N}{p_{1}}-1}_{2,p_{1},1})\cap \widetilde{L}^{1}(B^{\N+2,\frac{N}{p_{1}}+1}_{2,p_{1},1}) $ and $g^{n}\in \widetilde{L}^{\infty}(B^{\N,\frac{N}{p}}_{2,p,1})\cap \widetilde{L}^{1}(B^{\N+2,\frac{N}{p}+2}_{2,p,1}) $. We obtain then by proposition \ref{hybrid}:
$$
\begin{aligned}
 \|\big[\frac{q^{n}}{1+q^{n}}{\cal A}g^{n}\big]_{HF}\|_{\widetilde{L}^{1}(B^{\frac{N}{p}}_{p,1})}&\leq
\|T_{\frac{q^{n}}{1+q^{n}}}{\cal A}g^{n}\|_{\widetilde{L}^{1}(B^{\N,\frac{N}{p}}_{2,p,1})}+
\|T_{{\cal A}g^{n}}\frac{q^{n}}{1+q^{n}}\|_{\widetilde{L}^{1}(B^{\N-1,\frac{N}{p}}_{2,p,1})}\\
&\hspace{2cm}+\|R({\cal A}g^{n},\frac{q^{n}}{1+q^{n}})\|_{\widetilde{L}^{1}(B^{\N,\frac{N}{p}}_{2,p,1})},\\
& \leq C \|q^{n}\|_{\widetilde{L}^{\infty}(B^{\N-1,\frac{N}{p}}_{2,p,1})}\|{\cal A}g^{n}\|_{\widetilde{L}^{1}(B^{\N,\frac{N}{p}}_{2,p,1})},
\end{aligned}
$$
Next we have to use proposition \ref{hybrid} to treat the term $T_{{\cal A}h^{n}}(\frac{1}{1+q^{n}}-1)$ and $R({\cal A}h^{n},\frac{1}{1+q^{n}}-1)$ when $p_{1}>N$,
we have then:
$$\|T_{{\cal A}h^{n}}\frac{q^{n}}{1+q^{n}}\|_{\widetilde{L}^{1}(B^{\N,\frac{N}{p_{1}}-1}_{2,p_{1},1})}\leq
\|{\cal A}h^{n}\|_{\widetilde{L}^{1}(B^{\N,\frac{N}{p_{1}}-1}_{2,p_{1},1})}\|\frac{q^{n}}{1+q^{n}}\|_{\widetilde{L}^{\infty}(B^{\N-1,
\NN}_{2,p,1})},$$
where following the proposition \ref{hybrid}, we have chose $p=2$, $q=p_{1}$ , and as $p\geq p_{1}$ we have $\frac{1}{\lambda^{'}}=\frac{1}{p_{1}}-\frac{1}{p}$ and $\lambda=+\infty$. It means that: $\frac{N}{p_{1}}-1 \leq\NN$ (what is assumed) and $2\leq \lambda^{'}$ if $2\geq \frac{p_{1} p}{p-p_{1}}$. It means that we need of the following condition:
\begin{equation}
2\leq \frac{p_{1} p}{p-p_{1}}\;\;\;\mbox{and}\;\;\;\frac{N}{p_{1}}-1 \leq\NN.
 \label{condi1}
\end{equation}
Next we have as $\frac{N}{p_{1}}+\NN-1>0$ by proposition \ref{hybrid} for the rest term on the high frequencies:
$$\|\big[R({\cal A}h^{n},(\frac{1}{1+q^{n}}-1))\big]_{HF}\|_{\widetilde{L}^{1}(B^{\frac{N}{p_{1}}-1}_{p_{1},1})}\leq
\|h^{n}\|_{\widetilde{L}^{1}(B^{\N-1,\frac{N}{p_{1}}-1}_{2,p_{1},1})}\|\frac{1}{1+q^{n}}-1\|_{\widetilde{L}^{\infty}(B^{\N-1,
\NN}_{2,p,1})}.$$
We have seen that we need to treat this term of the condition:
\begin{equation}
\frac{N}{p_{1}}+\NN-1>0.
 \label{condi2}
\end{equation}
Easily we have by proposition \ref{hybrid} as $\widetilde{L}^{\infty}(B^{\N-1,
\NN}_{2,p,1})\h L^{\infty}$:
$$\|T_{\frac{q^{n}}{1+q^{n}}-1}{\cal A}h^{n}]_{HF}\|_{\widetilde{L}^{1}(B^{\frac{N}{p_{1}}-1}_{p_{1},1})}\leq
\|{\cal A}h^{n}\|_{\widetilde{L}^{1}(B^{\N-1,\frac{N}{p_{1}}-1}_{2,p_{1},1})}\|\frac{q^{n}}{1+q^{n}}\|_{\widetilde{L}^{\infty}(B^{\N-1,
\NN}_{2,p,1})}.$$
%$$
%\begin{aligned}
% \|\big[(\frac{1}{1+q^{n}}-1){\cal A}v_{1}^{n}\big]_{HF}\|_{\widetilde{L}^{1}(B^{\frac{N}{p_{1}}-1}_{p_{1},1})}&\leq
%\|(\frac{1}{1+q^{n}}-1){\cal A}v_{1}^{n}\|_{\widetilde{L}^{1}(B^{\N-1,\frac{N}{p_{1}}-1}_{2,p,1})},\\
%&\leq \|q^{n}\|_{\widetilde{L}^{\infty}(B^{\N-1,\frac{N}{p}}_{2,p,1})}\|{\cal %A}v_{1}^{n}\|_{\widetilde{L}^{1}(B^{\N-1,\frac{N}{p_{1}}-1}_{2,p_{1},1})}.
%\end{aligned}
%$$
We treat now the term $u^{n}\cdot\n u^{n}$ and we have as $u^{n}\in E$, it exists $h_{1}^{n}$ and $g^{n}_{1}$ such that $u^{n}=g^{n}_{1}+h^{n}_{1}$ with $h^{n}_{1}\in \widetilde{L}^{\infty}(\widetilde{B}^{\N-1,\NN}_{2,p,1})\cap \widetilde{L}^{1}(\widetilde{B}^{\N+1,\NN+1}_{2,p,1})$ and $g^{n}_{1}\in \widetilde{L}^{\infty}(\widetilde{B}^{\N-1,\frac{N}{p_{1}}-1}_{2,p_{1},1})\cap \widetilde{L}^{1}(\widetilde{B}^{\N+1,\NN+1}_{2,p,1})$.
We have then by proposition \ref{hybrid}:
$$
\begin{aligned}
&\|\big(h_{1}^{n}\cdot\n h_{1}^{n}\big)_{HF}\|_{\widetilde{L}^{1}(B^{\frac{N}{p}}_{p,1})}\leq
\|T_{h_{1}^{n}}\n h_{1}^{n}\|_{\widetilde{L}^{1}(\widetilde{B}^{\N,\NN}_{2,p,1})}
+\|T_{\n h_{1}^{n}}h_{1}^{n}\|_{\widetilde{L}^{1}(\widetilde{B}^{\N-1,\NN}_{2,p,1})}\\
&\hspace{8cm}+
\|R(\n h_{1}^{n},h_{1}^{n})\|_{\widetilde{L}^{1}(\widetilde{B}^{\N,\NN}_{2,p,1})},\\
&\hspace{4cm}\leq\| h_{1}^{n}\|_{\widetilde{L}^{1}(\widetilde{B}^{\N+1,\NN+1}_{2,p,1})}\| h_{1}^{n}\|_{\widetilde{L}^{\infty}(\widetilde{B}^{\N-1,\NN}_{2,p,1})}
\end{aligned}
$$
Next we have to treat the term $T_{g^{n}_{1}}\n g^{n}_{1}$ by using the proposition \ref{hybrid} with $\frac{1}{\lambda^{'}}=\frac{1}{p_{1}}-\frac{1}{p}$, $2\leq \lambda^{'}$ and $\frac{N}{p_{1}}-1\leq\NN$ then:
$$\|T_{g_{1}^{n}}\n g_{1}^{n}\|_{\widetilde{L}^{1}(\widetilde{B}^{\N-1,\frac{N}{p_{1}}-1}_{2,p_{1},1})}\leq
\|g_{1}^{n}\|_{\widetilde{L}^{\infty}(\widetilde{B}^{\N-1,\frac{N}{p_{1}}-1}_{2,p_{1},1})}
\|\n g_{1}^{n}\|_{\widetilde{L}^{1}(\widetilde{B}^{\N,\NN}_{2,p,1})}.$$
We have seen that we need of the conditions:
\begin{equation}
\frac{N}{p_{1}}-1\leq\NN\;\;\mbox{and}\;\;2\leq\frac{p_{1}p}{p-p_{1}}.
 \label{condi6}
\end{equation}
Easily we have as $\widetilde{B}^{\N,\NN}_{2,p,1}\h L^{\infty}$ by proposition \ref{hybrid}:
$$\|T_{\n g_{1}^{n}}g_{1}^{n}\|_{\widetilde{L}^{1}(\widetilde{B}^{\N-1,\frac{N}{p_{1}}-1}_{2,p_{1},1})}\leq
C\|g_{1}^{n}\|_{\widetilde{L}^{\infty}(\widetilde{B}^{\N-1,\frac{N}{p_{1}}-1}_{2,p_{1},1})}
\|\n g_{1}^{n}\|_{\widetilde{L}^{1}(\widetilde{B}^{\N,\NN}_{2,p,1})}.$$
To finish with the term $g^{n}_{1}\n g^{n}_{1}$, we have to treat the term $\big(R(g^{n}_{1},\n g^{n}_{1})\big)_{HF}$. By proposition \ref{hybrid}, as $\NN+\frac{N}{p_{1}}-1>0$ we have:
$$\|\big(R(g^{n}_{1},\n g^{n}_{1})\big)_{HF}\|_{\widetilde{L}^{1}(B^{\frac{N}{p_{1}}-1}_{p_{1},1})}\leq
C\|g_{1}^{n}\|_{\widetilde{L}^{\infty}(\widetilde{B}^{\N-1,\frac{N}{p_{1}}-1}_{2,p_{1},1})}
\|\n g_{1}^{n}\|_{\widetilde{L}^{1}(\widetilde{B}^{\N,\NN}_{2,p,1})}.$$
We have seen that we need of the conditions:
\begin{equation}
\frac{N}{p_{1}}-1+\NN> 0.
 \label{condi5}
\end{equation}
From the previous inequalities, we have obtained:
$$\|\big(g^{n}_{1}\cdot\n g^{n}_{1}\big)_{HF}\|_{\widetilde{L}^{1}(B^{\frac{N}{p_{1}}-1}_{p_{1},1})}\leq
C\|g_{1}^{n}\|_{\widetilde{L}^{\infty}(\widetilde{B}^{\N-1,\frac{N}{p_{1}}-1}_{2,p_{1},1})}
\|\n g_{1}^{n}\|_{\widetilde{L}^{1}(\widetilde{B}^{\N,\NN}_{2,p,1})}.$$
We can treat similarly the terms $g^{n}_{1}\cdot\n h^{n}_{1}$ and $h^{n}_{1}\cdot\n g^{n}_{1}$. We have finally under the conditions (\ref{condi1}), (\ref{condi2}), (\ref{condi6}) and (\ref{condi5}):
$$
\begin{aligned}
&\|\big(u^{n}\cdot\n u^{n}\big)_{HF}\|_{\widetilde{L}^{1}(B^{\frac{N}{p_{1}}-1}_{p_{1},1}+B^{\frac{N}{p}}_{p,1})}\leq
C\|u^{n}\|_{E}^{2}.
\end{aligned}
$$
We finish with the following term where $f$ is a regular function such that $f(0)=0$:
$$
\begin{aligned}
 &\|\big[\n(\D)^{-1}(P^{'}(\rho^{n}){\rm div}(\rho^{n}u^{n}))\big]_{HF}\|_{\widetilde{L}^{1}(B^{\frac{N}{p_{1}}-1}_{p_{1},1}+B^{\NN}_{p,1})}\\[2mm]
&\hspace{5cm}\leq
\|\big[\n(\D)^{-1}(f(q^{n}){\rm div}(q^{n}u^{n}))\big]_{HF}\|_{\widetilde{L}^{1}(B^{\frac{N}{p_{1}}-1}_{p_{1},1}+B^{\NN}_{p,1})}\\
&+
\|\big[\n(\D)^{-1}({\rm div}(q^{n}u^{n}))\big]_{HF}\|_{\widetilde{L}^{1}(B^{\frac{N}{p_{1}}-1}_{p_{1},1}+B^{\NN}_{p,1})}+\|\big[\n(\D)^{-1}{\rm div}(u^{n})\big]_{HF}\|_{\widetilde{L}^{1}(B^{\frac{N}{p_{1}}-1}_{p_{1},1}+B^{\NN}_{p,1})},\\
&\leq C \|{\rm div}(u^{n})\|_{\widetilde{L}^{1}(B^{\N,\NN+1}_{2,p,1})}(1+
\|q^{n}\|_{\widetilde{L}^{\infty}(B^{\N,\NN}_{2,p,1})}+
\|q^{n}\|^{2}_{\widetilde{L}^{\infty}(B^{\N,\NN}_{2,p,1})}).
\end{aligned}
$$
\\
\\
We have now to treat the case of low frequencies and in particular estimating $\|(F^{n})_{BF}\|_{\widetilde{L}^{1}(B^{\N-1}_{2,1})}
$ and $\|(G^{n})_{BF}\|_{\widetilde{L}^{1}(B^{\N-1}_{2,1})}$, we begin with $\|(F^{n})_{BF}\|_{\widetilde{L}^{1}(B^{\N-1}_{2,1})}$. We have then according proposition \ref{hybrid} if $p<\max(4,2N)$:
$$
\begin{aligned}
&\|(q^{n}{\rm div}u^{n})_{BF}\|_{\widetilde{L}^{1}(B^{\N-1}_{2,1})}\leq \|T_{q^{n}}({\rm div}u^{n})\|_{\widetilde{L}^{1}(\widetilde{B}^{\N-1,\NN-1}_{2,p,1})}+
\|T_{{\rm div}u^{n}}q^{n}\|_{\widetilde{L}^{1}(\widetilde{B}^{\N-1,\NN}_{2,p,1})}\\
&\hspace{8cm}+\|\big(R(q^{n},{\rm div}u^{n})\big)_{BF}\|_{\widetilde{L}^{1}(B^{\N-1}_{2,1})},\\
&\leq C\big(\|q^{n}\|_{\widetilde{L}^{2}(\widetilde{B}^{\N,\NN}_{2,p,1})}
\|u^{n}\|_{\widetilde{L}^{2}(\widetilde{B}^{\N,\frac{N}{p}}_{2,p,1})}+\|q^{n}\|_{\widetilde{L}^{\infty}(\widetilde{B}^{\N-1,\NN}_{2,p,1})}
\|u^{n}\|_{\widetilde{L}^{1}(\widetilde{B}^{\N+1,\frac{N}{p}+1}_{2,p,1})}\big).
\end{aligned}
$$
Here the only difficulty was to treat the term $R(q^{n},{\rm div}u^{n})$ when $N=2$, we need in this case of the previous condition:
\begin{equation}
p<\max(4,2N).
\label{condi4}
\end{equation}
and:
$$\|\big(R(q^{n},{\rm div}u^{n})\big)_{BF}\|_{\widetilde{L}^{1}(B^{\N-1}_{2,1})}\leq
C\|q^{n}\|_{\widetilde{L}^{2}(\widetilde{B}^{\N,\NN}_{2,p,1})}
\|u^{n}\|_{\widetilde{L}^{2}(\widetilde{B}^{\N,\frac{N}{p}}_{2,p,1})}.$$
Similarly we have by using proposition \ref{hybrid} with condition (\ref{condi4}):
$$\|(u^{n}\cdot\n q^{n})_{BF}\|_{\widetilde{L}^{1}(B^{\N-1}_{2,1})}\leq
C\|q^{n}\|_{\widetilde{L}^{2}(\widetilde{B}^{\N,\NN}_{2,p,1})}
\|u^{n}\|_{\widetilde{L}^{\infty}(\widetilde{B}^{\N,\frac{N}{p}}_{2,p,1})}.$$
%$$
%\begin{aligned}
%\|q^{n}{\rm div}u^{n}\|_{\widetilde{L}^{1}(\widetilde{B}^{\N-1,\frac{N}{p}-1}_{2,p,1})}&\leq C\|q^{n}\|_{\widetilde{L}^{\infty}(\widetilde{B}^{\N-1,\frac{N}{p_{1}}-2}_{2,p_{1},1})}\|u^{n}\|_{\widetilde{L}^{1}(\widetilde{B}^{\N+1,\frac{N}{p_{1}}+1}_{2,p_{1},1}+\widetilde{B}^{\N+1,\frac{N}{p}+1}_{2,p,1})},
%\end{aligned}
%$
%As $g(q^{n})=P(1+q^{n})-P(1)$, by lemma \ref{composition} we have:
%$$
%\begin{aligned}
%\|P^{'}(1)q^{n}-(1+q^{n})g(q^{n})\|_{L^{1}_{T}(B^{\frac{N}{p}-1,\frac{N}{p}}_{p,1})}&\leq\|P^{'}(1)(q^{n}\times %q^{n})\|_{L^{1}_{T}(B^{\frac{N}{p}-1,\frac{N}{p}}_{p,1})}\\
%&\leq  C\|q^{n}\|_{L^{1}_{T}(B^{\frac{N}{p}-1,\frac{N}{p}}_{p,1})}\|q^{n}\|_{\widetilde{L}^{\infty}(B^{\NN}_{p,1})}.
%\end{aligned}
%$$
\\
\\
We now want to estimate $\|G^{n}\|_{\widetilde{L}^{1}(\widetilde{B}^{\N-1,\frac{N}{p_{1}}-1}_{2,p_{1},1}+\widetilde{B}^{\N-1,\frac{N}{p}}_{2,p,1})}$, we begin with $\|(\frac{q^{n}}{1+q^{n}}
{\cal A}u^{n})_{HF}\|_{\widetilde{L}^{1}(B^{\N-1}_{2,1})}$, the main difficulty corresponds to treat $T_{{\cal A}u^{n}}\frac{q^{n}}{1+q^{n}}$ and $R(\frac{q^{n}}{1+q^{n}},{\cal A}u^{n})$. We have by using proposition \ref{hybrid} if:
$\frac{1}{2}\leq\frac{2}{p}$, $N-1>0$, $2\NN-1>0$. We recall here that $\widetilde{L}^{\infty}(\widetilde{B}^{\N,\NN}_{2,p,1})\h \widetilde{L}^{\infty}(\widetilde{B}^{\N-1,\NN}_{2,p,1})$, we have then:
$$
\begin{aligned}
&\|\big(R(\frac{q^{n}}{1+q^{n}},
{\cal A}u^{n})\big)_{BF}\|_{\widetilde{L}^{1}(B^{\N-1}_{2,1})}\leq C\|\frac{q^{n}}{1+q^{n}}\|_{\widetilde{L}^{\infty}(\widetilde{B}^{\N,\NN}_{2,p,1})}\|{\cal A}u^{n}\|_{\widetilde{L}^{1}(\widetilde{B}^{\N-1,\frac{N}{p}-1}_{2,p,1})}.
\end{aligned}
$$
We need then of the following conditions:
\begin{equation}
p< \max(4,2N).
 \label{condi3}
\end{equation}
%$$
%\begin{aligned}
%&\|\frac{q^{n}}{1+q^{n}}
%{\cal A}u^{n}\|_{\widetilde{L}^{1}(\widetilde{B}^{\N-1,\frac{N}{p_{1}}-1}_{2,p_{1},1}+\widetilde{B}^{\N-1,\frac{N}{p}}_{2,p,1})}
%\leq C\|\frac{q^{n}}{1+q^{n}}\|_{\widetilde{L}^{\infty}(\widetilde{B}^{\N-1,\NN}_{2,p,1})}\big(\|d_{1}^{n}\|_{\widetilde{L}^{1}(\widetilde{B}^{\N+1,\frac{N}{p_{1}}+1}_{2,p_{1},1}+\widetilde{B}^{\N+1,\frac{N}{p}+2}_{2,p,1})}\\
%&+\|(q^{n})_{HF}\|_{\widetilde{L}^{1}(B^{\NN}_{p,1})}+
%\|d_{1}^{n}\|_{\widetilde{L}^{1}(\widetilde{B}^{\N+1,\frac{N}{p_{1}}+1}_{2,p_{1},1}+\widetilde{B}^{\N+1,\frac{N}{p}+2}_{2,p,1})}\big)
%\end{aligned}
%$$
%Next we have:
%$$\|\n(\D)^{-1}(P^{'}(\rho^{n}){\rm div}(\rho^{n}u^{n}))\|_{\widetilde{L}^{1}(B^{\NN}_{p,1}+B^{\frac{N}{p_{1}}-1}_{p_{1},1})}$$
%Next we have:
%$$
%\begin{aligned}
%\|u^{n}\cdot\n u^{n}\|_{\widetilde{L}^{1}(B^{\NN}_{p,1}+B^{\frac{N}{p_{1}}-1}_{p_{1},1})}\leq
%\end{aligned}
%$$
%Moreover we recall that according to proposition \ref{important}:
%$$\|q^{n}\|_{L^{2}_{T}(B^{\frac{N}{2}})}^{2}\leq\|q^{n}\|_{L_{T}^{\infty}(\widetilde{B}^{\frac{N}{2}-1,\frac{N}{2}})}
%\|q^{n}\|_{L_{T}^{1}(\widetilde{B}^{\frac{N}{2}+1,\frac{N}{2}})}.
%$$
%We proceed similarly to estimate
%$\|H_{1}^{n}\|_{L_{T}^{1}(B^{\frac{N}{2}-1})}$ and
Next we have according proposition \ref{hybrid} with $\lambda=\lambda^{'}=+\infty$:
$$
\begin{aligned}
\|\big(T_{{\cal A}u^{n}}\frac{q^{n}}{1+q^{n}}\big)_{BF}\|_{\widetilde{L}^{1}(B^{\N-1}_{2,1})}&\leq
\|T_{{\cal A}u^{n}}\frac{q^{n}}{1+q^{n}}\|_{\widetilde{L}^{1}(\widetilde{B}^{\N-1,\NN-1}_{2,p,1})},\\
&\leq C\|\frac{q^{n}}{1+q^{n}}\|_{\widetilde{L}^{\infty}(\widetilde{B}^{\N,\NN}_{2,p,1})}\|{\cal A}u^{n}\|_{\widetilde{L}^{1}(\widetilde{B}^{\N-1,\frac{N}{p}-1}_{2,p,1})}
\end{aligned}
$$
Finally we have:
$$\|(K(q^{n})\n q^{n})_{BF}\|_{\widetilde{L}^{1}(B^{\N-1}_{2,1})}\leq \|q^{n}\|_{\widetilde{L}^{2}(\widetilde{B}^{\N,\NN}_{2,p,1})}.$$
To finish, it stays in low frequencies the terms: $T_{\n u^{n}}u^{n}$, $T_{u^{n}}\n u^{n}$ and $R(u^{n},\n u^{n})$. We have then by proposition \ref{hybrid} if $p<\max(4,2N))$:
$$\|(R(h_{1}^{n},\n h_{1}^{n}))_{BF}\|_{\widetilde{L}^{1}(B^{\N-1}_{2,1})}\leq C\|h_{1}^{n}\|_{\widetilde{L}^{\infty}(\widetilde{B}^{\N-1,\NN-1}_{2,p,1})}\|\n h_{1}^{n}\|_{\widetilde{L}^{1}(\widetilde{B}^{\N,\NN}_{2,p,1})}.$$
We have seen that we need again of condition (\ref{condi4}).\\
We want treat now $(R(g_{1}^{n},\n g_{1}^{n}))_{BF}$, we have then by proposition \ref{hybrid} if $\frac{N}{p_{1}}+\NN-1>0$ and
$\frac{1}{2}\leq\frac{1}{p}+\frac{1}{p_{1}}$.
$$\|(R(g_{1}^{n},\n g_{1}^{n}))_{BF}\|_{\widetilde{L}^{1}(B^{\N-1}_{2,1})}\leq C\|g_{1}^{n}\|_{\widetilde{L}^{\infty}(\widetilde{B}^{\N-1,\frac{N}{p_{1}}-1}_{2,p-{1},1})}\|\n g_{1}^{n}\|_{\widetilde{L}^{1}(\widetilde{B}^{\N,\NN}_{2,p,1})}.$$
We have seen that we need of the following conditions:
\begin{equation}
\frac{1}{2}\leq\frac{1}{p}+\frac{1}{p_{1}}\;\;\mbox{and}\;\;\frac{N}{p_{1}}+\NN-1>0.
 \label{condi9}
\end{equation}
Next we have by proposition \ref{hybrid} if:
$$\|T_{\n h_{1}^{n}}h_{1}^{n}\|_{\widetilde{L}^{1}(\widetilde{B}^{\N-1,\NN}_{2,p,1})}\leq C\|h_{1}^{n}\|_{\widetilde{L}^{\infty}(\widetilde{B}^{\N-1,\NN}_{2,p,1})}\|\n h_{1}^{n}\|_{\widetilde{L}^{1}(\widetilde{B}^{\N,\NN}_{2,p,1})}.$$
and
$$\|T_{\n g_{1}^{n}}g_{1}^{n}\|_{\widetilde{L}^{1}(\widetilde{B}^{\N-1,\NN}_{2,p,1})}\leq C\|g_{1}^{n}\|_{\widetilde{L}^{\infty}(\widetilde{B}^{\N-1,\frac{N}{p_{1}}-1}_{2,p-{1},1})}\|\n g_{1}^{n}\|_{\widetilde{L}^{1}(\widetilde{B}^{\N,\NN}_{2,p,1})}.$$
We proced similarly to control $T_{u^{n}}\n u^{n}$.
Therefore the above inequalities with conditions (\ref{condi1}), (\ref{condi2}), (\ref{condi6}), (\ref{condi5}), (\ref{condi4}), (\ref{condi3}) and (\ref{condi9}) imply that for all $t\in\R$ we have :
$$\|(q^{n},u^{n})\|_{H_{t}}\leq Ce^{C\|(q^{n},u^{n})\|_{H_{t}}}(\|q_{0}\|_{B^{\N-1\NN}_{2,p,1}}+\|u_{0}\|_{B^{\N-1,\frac{N}{p_{1}}-1}_{2,p_{1},1}}+\|f\|_{\widetilde{L}^{1}(B^{\N-1\frac{N}{p_{1}}-1
}_{2,p_{1},1}}+
\|(q^{n},u^{n})\|^{2}_{H_{t}}).$$
From a standard bootstrap argument, it is now easy to conclude that there exists a positive constant $c$ such that if the data has been chosen so small as to satisfy:
$$\|q_{0}\|_{B^{\N-1\NN}_{2,p,1}}+\|u_{0}\|_{B^{\N-1,\frac{N}{p_{1}}-1}_{2,p_{1},1}}+\|f\|_{\widetilde{L}^{1}(B^{\N-1\frac{N}{p_{1}}-1
}_{2,p_{1},1}}\leq c.$$
then it exists $C>0$ such that for all $t\in \R$:
$$\|(q^{n},u^{n})\|_{H_{t}}\leq C,\;\;\forall t\in\R.$$
\subsubsection*{Compactness arguments}
Let us first focus on the convergence of $(q^{n})_{n\in\mathbb{N}}$. We claim that, up to extraction, $(q^{n})_{n\in\mathbb{N}}$ converges in the distributional sense to some function $q$ such that:
\begin{equation}
 q\in\widetilde{L}^{\infty}(B^{\N-1,\NN}_{2,p,1})\cap \widetilde{L}^{1}(B^{\N+1,\NN}_{2,p,1}).
\label{49}
\end{equation}
%We now have to prove that $(a^{n},u^{n})_{n\in\mathbb{N}}$ tends (up to a subsequence) to some function $(a,u)$ which belongs to $E_{T}$. Here we recall that:
%$$E_{T}=\widetilde{C}([0,T],B^{\NN}_{p,1})\times\big(\widetilde{L}^{\infty}(B^{\frac{N}{p_{1}}-1}_{p_{1},1}+B^{\NN+1}_{p,1})\cap \widetilde{L}^{1}(
%B^{\frac{N}{p_{1}}+1}_{p_{1},1}+B^{\NN+2}_{p,1})\big).$$
%As $u^{n}=v_{1}^{n}+\frac{1}{\nu}v^{n}$ with ${\rm div}v^{n}=P(\rho^{n})$ it would be esay to conclude that up to a subsequence
% $(u^{n})_{n\in\mathbb{N}}$ tends to some function $u$ where $(a,u)$ and satisfies (\ref{0.6}).
The proof is based on Ascoli's theorem and compact embedding for Besov spaces. As similar arguments have been employed in \cite{DL} or \cite{DW}, we only give the outlines of the proof.
%\begin{itemize}
%\item Convergence of $(a^{n})_{n\in\mathbb{N}}$:\\
We may write that:
%We use the fact that $\widetilde{a}^{n}=a^{n}-a^{n}_{0}$ satisfies:
$$\p_{t} q^{n}=-u^{n}\cdot\n q^{n}-(1+q^{n}){\rm div}u^{n}.$$
Since $(u^{n})_{n\in\mathbb{N}}$ is uniformly bounded in $\widetilde{L}^{2}_{T}(B^{\N,\NN+1}_{2,p,1}+B^{\N,\frac{N}{p_{1}}}_{2,p_{1},1})$ and $q^{n}\in \widetilde{L}^{\infty}(B^{\N-1,\NN}_{2,p,1})$, we have $(1+q^{n}){\rm div}u^{n}$ which is bounded in
$\widetilde{L}^{2}_{T}(B^{\N-1,\NN}_{2,p,1}+B^{\N-1,\frac{N}{p_{1}}-1}_{2,p_{1},1})$ with the conditions between $p$ and $p_{1}$ in theorem \ref{theo1}. Similarly $u^{n}\cdot\n q^{n}$ is bounded in $\widetilde{L}^{2}_{T}(B^{\N-1,\NN}_{2,p,1}+B^{\N-1,\frac{N}{p_{1}}-1}_{2,p_{1},1})$. Finally as $p\geq p_{1}$, we have proved that $\p_{t}q^{n}$ is bounded in $\widetilde{L}^{2}_{T}(B^{\N-1,\NN-1}_{2,p,1})$, it means that $(q^{n})_{n\in\mathbb{N}}$  seen as a sequence of $B^{\N-1,\NN-1}_{2,p,1}$ valued functions is equicontinuous in $\R$.
In addition $(q^{n})_{n\in\mathbb{N}}$ is bounded in $C(\R,B^{\N-1,\NN-1}_{2,p,1}\cap B^{\N,\NN}_{2,p,1})$. As the embedding $B^{\N-1,\NN-1}_{2,p,1}\cap B^{\N,\NN}_{2,p,1}$ is locally compact (see \cite{BCD}, Chap2), one can thus conclude by means of Ascoli's theorem and Cantor diagonal extraction process that there exists some distribution $q$ such that up to an omitted extraction $(\psi q^{n})_{n\in\mathbb{N}}$ converges to $\psi q$
in $\C(\R,B^{\N-1,\NN-1}_{2,p,1})$ for all smooth $\psi$ with compact support in $\R^{=}\times\R^{N}$. Then by using the so-called Fatou property for the Besov spaces, one can conclude that (\ref{49}) is satisfied. (the reader may consult \cite{BCD}, Chap 10 too).
By proceeding similarly, we can prove that up to extraction, $(u^{n})_{n\in\mathbb{N}}$ converges in the distributional sense to some function $u$ such that:
\begin{equation}
 u\in\widetilde{L}^{\infty}(B^{\N-1,\frac{N}{p_{1}}-1}_{2,p_{1},1}+B^{\N-1,\NN}_{2,p_{1},1})\cap \widetilde{L}^{1}(B^{\N+1,\NN+1}_{2,p,1}).
\label{51}
\end{equation}
In order to complete the proof of the existence part of theorem \ref{theo1}, it is only a matter of checking the continuity properties with respect to time, namely that:
$$
\begin{aligned}
&q\in\widetilde{C}(\R^{+},\widetilde{B}^{\N-1,\NN}_{2,p,1})\;\;\mbox{and}\;\;u\in \widetilde{C}(\R^{+},\widetilde{B}^{\N-1,\NN}_{2,p,1}+\widetilde{B}^{\N-1,\frac{N}{p_{1}}-1}_{2,p_{1},1}).
\end{aligned}
$$
As regards $q$, it suffices to notice that, according to (\ref{49}), (\ref{51}) and to the product laws in the Besov spaces, we have:
$$\p_{t}q+u\cdot q=-(1+q){\rm div}u\in \widetilde{L}^{1}(B^{\N,NN}_{2,p,1}.$$
As $q_{0}\in B^{\N,\NN}_{2,p,1}$, classical results for the transport equation (see \cite{BCD}, Chap 3) ensure that $q\in
\widetilde{C}(\R^{+},\widetilde{B}^{\N,\NN}_{2,p,1})$. And as previously, we have shown that $q\in
\widetilde{C}(\R^{+},\widetilde{B}^{\N-1,\NN-1}_{2,p,1})$, it means clearly that $q\in
\widetilde{C}(\R^{+},\widetilde{B}^{\N-1,\NN}_{2,p,1})$.\\
For getting the continuity result for $u$, one may similarly use the properties of the heat equation on $v_{1}$ in high frequncies and on $u$ in low frequencies.
\subsection*{The proof of the uniqueness}
In the case $\frac{2}{N}\leq\frac{1}{p}+\frac{1}{p_{1}}$, the uniqueness has been established in \cite{DL,H2}.
%\subsection{Proof of corollary \ref{corollary1}}
%We have just to remark than $L^{2}$ is include in low frequencies in $B^{\N-1}_{2,1}$ and than the Orlicz space $L^{\gamma}_{2}$ is include in $B^{\N-1}_{2,1}$ in low frequencies too.
\section{Proof of corollary \ref{corollaire1} and \ref{corollaire2}}
\label{section5}
\subsection{Proof of corollary \ref{corollaire1}}
We want here to avoid the condition $p<\max(4,2N)$. For simplicity we will treat only the case $N=3$. This condition appears when we want treat the terms of rest in low frequencies. For
resolving this problem as in the paper of F. Charve and R. Danchin in \cite{CD}, we need of additionnal condition in high frequencies on $q_{0}$ and $u_{0}$.\\
We want then to follow the same strategy as in the proof of theorem \ref{theo1}. It means that we use the same
standard scheme which consists in the construction of approximate solutions, some uniform bounds and results of compactness. We will use the same notations as in proofof theorem  \ref{theo1}. We just want treat the non linear term where appears the condition $p<\max(2N,4)$ in an other way by using the additional hypothesis that we have on $(q_{0},u_{0})\in \widetilde{B}^{0,1}_{2,\infty}\times B^{0}_{2,\infty}$. The rest of the proof will be the same as in theorem \ref{theo1}. We will work with the same space as in the proof \ref{theo1} except that we attend additional regularity on $(q^{n},u^{n})$ in $E^{'}$ with:
$$E_{1}^{'}=(\widetilde{L}^{\infty}(\widetilde{B}^{0,1}_{2,\infty}\cap \widetilde{L}^{1}(\widetilde{B}^{2,1}_{2,\infty})\times
(\widetilde{L}^{\infty}(\widetilde{B}^{0}_{2,\infty}\cap \widetilde{L}^{1}(\widetilde{B}^{2}_{2,\infty}).$$
Here $(q^{n},u^{n})_{n\in\mathbb{N}}$ is the solution of the following system:
\begin{equation}
\begin{cases}
\begin{aligned}
&\p_{t}q^{n}+u^{n}\cdot\n q^{n}+{\rm div}u^{n}=F^{n}\\
&\p_{t}u^{n}+u^{n}\cdot\n u^{n}-{\cal A}u^{n}+P^{'}(1)\n q^{n}=G^{n}+f^{n}\\
&(q^{n},u^{n})_{/t=0}=(q^{n}_{0},u^{n}_{0}),
\end{aligned}
\end{cases}
\label{5grandeequation}
\end{equation}
which verifies proposition 4 in \cite{CD} with:
$$
\begin{aligned}
&F^{n}=-q^{n}{\rm div}u^{n},\\
&G^{n}=-\frac{q^{n}}{1+q^{n}}{\cal A}u^{n}+(P^{'}(1)-P^{'}(1+q^{n}))\n q^{n}.
\end{aligned}
$$
%Let us set:
%$$
%\begin{aligned}
%&\beta=\|u_{0}\|_{\widetilde{B}^{\N-1,\frac{N}{p_{1}}-1}_{2,p_{1},1}}+
%\|q_{0}\|_{\widetilde{B}^{\N-1,\frac{N}{p}}_{2,p,1}}+
%\|f\|_{L^{1}_{T}(B^{\frac{N}{p_{1}}-1}_{p_{1},1})},\\
%&E(q,u)=\|q\|_{\widetilde{L}^{\infty}(\widetilde{B}^{\N-1,\NN}_{2,p,1})}+\|u\|_{\widetilde{L}^{\infty}(\widetilde{B}^{\N-1, \frac{N}{p_{1}}-1}_{2,p_{1},1}+\widetilde{B}^{\N-1, \frac{N}{p}}_{2,p,1})}+\|q\|_{\widetilde{L}^{1}(\widetilde{B}^{\N+1,\NN}_{2,p,1})}\\
%&\hspace{8cm}+
%\|u\|_{\widetilde{L}^{1}(\widetilde{B}^{\N+1,\frac{N}{p_{1}}+1}_{2,p_{1},1}+\widetilde{B}^{\N+1,\frac{N}{p}+1}_{2,p,1})},\\
%&E_{1}(q,u)=\|q\|_{\widetilde{L}^{\infty}(B^{\N-1}_{2,1})}+\|u\|_{\widetilde{L}^{\infty}(B^{\N-1}_{2,1})}+
%\|q\|_{\widetilde{L}^{1}(B^{\N+1}_{2,1})}+
%\|u\|_{\widetilde{L}^{1}(B^{\N+1}_{2,1})}.\\
%&E_{2}(q,u)=\|q\|_{\widetilde{L}^{\infty}(B^{\NN}_{p,1})}+\|u\|_{\widetilde{L}^{\infty}(B^{ \frac{N}{p_{1}}-1}_{p_{1},1}+B^{ \frac{N}{p}}_{p,1})}+\|q\|_{\widetilde{L}^{1}(B^{\NN}_{p,1})}\\
%&\hspace{8cm}+
%\|u\|_{\widetilde{L}^{1}(B^{\frac{N}{p_{1}}+1}_{p_{1},1}+B^{\frac{N}{p}+2}_{p,1})}.
%\|q\|_{\widetilde{L}_{t}^{\infty}(\widetilde{B}^{\NN-1,\NN}_{p,1})}+\|q\|_{\widetilde{L}_{t}^{1}(\widetilde{B}^{\NN-1,\NN}_{p,1})}
%+\|v_{1}\|_{\widetilde{L}_{t}^{\infty}(B^{\NN-1}_{p,1}+B^{\frac{N}{p_{1}}-1}_{p_{1},1})}\\
%&\hspace{10cm}+\|v_{1}\|_{\widetilde{L}^{1}_{t}(
%B^{\NN+1}_{p,1}+B^{\frac{N}{p_{1}}+1}_{p_{1},1})}.
%\end{aligned}
%$$
We apply exactly the same proof than for theorem \ref{theo1}, however we have to complete the uniform bounds by showing that $(q^{n},u^{n})$ is uniformly bounded in $H^{'}\cap E^{'}_{1}$, moreover we have to treat differently the term in low frequencies where appears the conditions $p<\max(4,2N)$ and $\frac{1}{2}\leq\frac{1}{p}+\frac{1}{p_{1}}$ by using the fact that $(q^{n},u^{n})$ in $E^{'}_{1}$.\\
We begin with treatinf the terms $\|F^{n}\|_{\widetilde{L}^{1}(B^{0,1}_{2,\infty})}$ and
$\|G^{n}\|_{\widetilde{L}^{1}(B^{0}_{2,\infty})}$
 by using properties
of continuity on the paraproduct and proposition 4 of \cite{VD}. We have then:
$$\|T^{'}_{q^{n}}{\rm div}u^{n}\|_{\widetilde{L}^{1}(B^{1}_{2,\infty})}\leq
\|q^{n}\|_{L^{\infty}(L^{\infty})}\|{\rm div}u^{n}\|_{\widetilde{L}^{1}(B^{1}_{2,\infty})}.$$
Similarly:
$$\|T^{'}_{q^{n}}{\rm div}u^{n}\|_{\widetilde{L}^{1}(B^{0}_{2,\infty})}\leq \|q^{n}\|_{L^{2}(L^{\infty})}\|{\rm div}u^{n}\|_{\widetilde{L}^{2}(B^{0}_{2,\infty})}.
$$
Next we have:
$$\|K(q)\n q\|_{\widetilde{L}^{1}(B^{0}_{2,\infty})}\leq\|q\|_{L^{2}(L^{\infty})}\|K(q)\|_{\widetilde{L}^{2}(B^{1}_{2,\infty})}.$$
For the term $(\frac{1}{\rho}-\frac{1}{\bar{\rho}})\D u=J(q)\D u$ with $J$ regular and $J(0)=0$, we have:
$$
\begin{aligned}
&\|T_{J(q)}\D u\|_{\widetilde{L}^{1}(B^{0}_{2,\infty})}\leq\|J(q)\|_{L^{\infty}(L^{\infty})}\|\D u\|_{\widetilde{L}^{1}(B^{0}_{2,\infty})},\\
&\|T_{\D u}J(q)\|_{\widetilde{L}^{1}(B^{0}_{2,\infty})}\leq\|J(q)\|_{\widetilde{L}^{\infty}(B^{1}_{2,\infty})}\|\D u\|_{\widetilde{L}^{1}(\widetilde{B}^{\N-1,\NN-1}_{2,p,\infty})},
\end{aligned}
$$
Concerning the remainder, we have if $p\geq 2$:
$$\|R(J(q),\D u)\|_{\widetilde{L}^{1}(B^{0}_{2,\infty})}\leq \|J(q)\|_{\widetilde{L}^{\infty}(B^{1}_{2,\infty})}\|\D u\|_{\widetilde{L}^{1}(\widetilde{B}^{\N-1,\NN-1}_{2,p,\infty})}.$$
We have then obtained:
$$\|J(q)\D u\|_{\widetilde{L}^{1}(B^{0}_{2,\infty})}\leq \|q\|_{\widetilde{L}^{\infty}(B^{1}_{2,\infty})}\| u\|_{\widetilde{L}^{1}(\widetilde{B}^{\N+1,\NN+1}_{2,p,\infty})}+\|J(q)\|_{L^{\infty}(L^{\infty})}\|\D u\|_{\widetilde{L}^{1}(B^{0}_{2,\infty})}.$$
It stays now to treat now the terms where appears the conditions $p<\max(4,2N)$ and $\frac{1}{2}\leq\frac{1}{p}+\frac{1}{p_{1}}$ in an other way.
%$$
%\begin{aligned}
%\|\big(R({\cal A}u^{n},\frac{q^{n}}{1+q^{n}})\big)_{BF}\|_{\widetilde{L}^{1}(B^{0}_{2,1})}&\leq C\|\big(R({\cal A}u^{n},\frac{q^{n}}{1+q^{n}})\big)_{BF}\|_{\widetilde{L}^{1}(B^{\N-1}_{2,1})},\\
%&\leq C\|q^{n}\|_{\widetilde{L}^{\infty}(\widetilde{B}^{1}_{2,2})}
%\|u^{n}\|_{\widetilde{L}^{1}(\widetilde{B}^{\N+1,\frac{N}{p}+1}_{2,p,1})}.
%\end{aligned}
%$$
%Next we have as $\widetilde{B}^{\N,\frac{N}{p}+1}_{2,p,1})+\widetilde{B}^{\N,\frac{N}{p_{1}}}_{2,p_{1},1})\h \widetilde{B}^{\N,\frac{N}{p}}_{2,p,1})$:
As $u^{n}\in \widetilde{L}^{2}(B^{\N,\NN+1}_{2,p,1}+B^{\N,\frac{N}{p_{1}}}_{2,p_{1},1})$, we set $u^{n}=g^{n}+h^{n}$ with $g^{n}\in
\widetilde{L}^{2}(B^{\N,\NN+1}_{2,p,1})$ and $h^{n}\in \widetilde{L}^{2}(B^{\N,\frac{N}{p_{1}}}_{2,p_{1},1})$. According to proposition \ref{hybrid}, we have:
$$
\begin{aligned}
\|\big(R(q^{n},{\rm div}g^{n})\big)_{BF}\|_{\widetilde{L}^{1}(B^{\N-1,}_{2,1})}&\leq C \|\big(R(q^{n},{\rm div}g^{n})\big)_{BF}\|_{\widetilde{L}^{1}(B^{0}_{2,1})},\\
&\leq C\|q^{n}\|_{\widetilde{L}^{2}(\widetilde{B}^{1}_{2,\infty})}
\|{\rm div}g^{n}\|_{\widetilde{L}^{2}(\widetilde{B}^{\N-1,\frac{N}{p}-1}_{2,p,1})},\\
&\leq C\|q^{n}\|_{\widetilde{L}^{2}(\widetilde{B}^{1}_{2,\infty})}
\|{\rm div}g^{n}\|_{\widetilde{L}^{2}(\widetilde{B}^{\N-1,\frac{N}{p}}_{2,p,1})}.
\end{aligned}
$$
$$
\begin{aligned}
\|\big(R(q^{n},{\rm div}h^{n})\big)_{BF}\|_{\widetilde{L}^{1}(B^{\N-1,}_{2,1})}&\leq C \|\big(R(q^{n},{\rm div}h^{n})\big)_{BF}\|_{\widetilde{L}^{1}(B^{0}_{2,1})},\\
&\leq C\|q^{n}\|_{\widetilde{L}^{2}(\widetilde{B}^{1}_{2,\infty})}
\|{\rm div}h^{n}\|_{\widetilde{L}^{2}(\widetilde{B}^{\N-1,\frac{N}{p_{1}}-1}_{2,p_{1},1})}.
\end{aligned}
$$
Next we have:
$$
\begin{aligned}
\|\big(R(u^{n},\n u^{n})\big)_{BF}\|_{\widetilde{L}^{1}(B^{\N-1,}_{2,1})}&\leq C \|\big(R(u^{n},\n u^{n})\big)_{BF}\|_{\widetilde{L}^{1}(L^{B^{0}_{2,1}})},\\
&\leq C
\|u^{n}\|_{\widetilde{L}^{\infty}(B^{0}_{2,\infty})}
\|\n u^{n}\|_{\widetilde{L}^{1}(\widetilde{B}^{\N,\frac{N}{p}}_{2,p,1})}.
\end{aligned}
$$
It stays to control the term $R(\frac{q^{n}}{1+q^{n}},\D u^{n})$ in low frequencies:
$$
\begin{aligned}
\|\big(R(\frac{q^{n}}{1+q^{n}},\D u^{n})\big)_{BF}\|_{\widetilde{L}^{1}(B^{\N-1,}_{2,1})}&\leq C \|\big(R(\frac{q^{n}}{1+q^{n}},\D u^{n})\big)_{BF}\|_{\widetilde{L}^{1}(L^{2})},\\
&\leq C
\|q^{n}\|_{\widetilde{L}^{\infty}(B^{1}_{2,\infty})}
\|\D u^{n}\|_{\widetilde{L}^{1}(\widetilde{B}^{\N-1,\frac{N}{p}-1
}_{2,p,1})}.
\end{aligned}
$$
Therefore the above nequalities imply that for all $t\in[0,T]$ we have :
$$
\begin{aligned}
&\|(q^{n},u^{n})\|_{H_{t}\cap (E_{1}^{'})_{t}}\leq Ce^{C\|(q^{n},u^{n})\|_{H_{t}\cap (E_{1}^{'})_{t}}}(\|q_{0}\|_{B^{\N-1\NN}_{2,p,1}\cap B^{0,1}_{2,r}}+\|u_{0}\|_{B^{\N-1,\frac{N}{p_{1}}-1}_{2,p_{1},1}\cap B^{0}_{2,r}}\\
&\hspace{5cm}+\|f\|_{\widetilde{L}^{1}(B^{\N-1,\frac{N}{p_{1}}-1}_{2,p_{1},1}\cap B^{0}_{2,r}}+
\|(q^{n},u^{n})\|^{2}_{H_{t}\cap (E_{1}^{'})_{t}}).
\end{aligned}
$$
From a standard bootstrap argument, it is now easy to conclude that there exists a positive constant $c$ such that if the data has been chosen so small as to satisfy:
$$\|q_{0}\|_{B^{\N-1\NN}_{2,p,1}\cap B^{0,1}_{2,r}}+\|u_{0}\|_{B^{\N-1,\frac{N}{p_{1}}-1}_{2,p_{1},1}\cap B^{0}_{2,r}}+\|f\|_{\widetilde{L}^{1}(B^{\N-1,\frac{N}{p_{1}}-1}_{2,p_{1},1}\cap B^{0}_{2,r}}\leq c.$$
then it exists $C>0$ such that for all $t\in \R$:
$$\|(q^{n},u^{n})\|_{H_{t}\cap E^{'}_{t}}\leq C,\;\;\forall t\in\R.$$
To conclude we follow the previous proof of theorem \ref{theo1}. Compactness results go along the lines
of the proof of theorem \ref{theo1}.
\subsection{Proof of corollary \ref{corollaire2}}
We follow here exactly the lines
of the proof of theorem \ref{theo1} except that we introduce a new effective velocity. Indeed in our case $v$ verifies the following elliptic equation:
$${\rm}(\mu(1)D v)+\n(\lambda(1){\rm div}v)=\n P(\rho)+{\rm}(f_{1}(q) D v)+\n(f_{2}(q){\rm div}v),$$
with $f_{1}(q)=\mu(1+q)-\mu(1)$ and $f_{2}(q)=\lambda(1)-\lambda(1+q)$. We can resolve this elliptic equation as $\mu(1)\geq c>0$ and $\mu(1)+\lambda(1)\geq c>0$, indeed in our case we work away from the vacuum.
To do this we have to use the estimates on the Lam\'e operator of the appendix in \cite{H}. More precisely we have as $q\in\widetilde{L}^{\infty}(B^{\N-1,\NN}_{2,p,1})$ for $r\geq 1$, $p,q\geq1$ and $|s_{1}|<\N$, $|s_{2}|<\NN$<:
$$\|v\|_{\widetilde{L}^{r}(B^{s_{1},s_{2}}_{p,q,1})}\leq C\|q\|_{\widetilde{L}^{r}(B^{s_{1}-1,s_{2}-1}_{p,q,1})}.$$
Indeed as $q$ is small, the terms of rest with $f_{1}(q)$ and $f_{2}(q)$ are easy to treat.
It means as in the proof of theorem \ref{theo1}, $v$ is one derivative more regular than $q$ in high frequencies and that we can estimate $v$ in function of $q$.
Moreover we have $\p_{t}v$ which verifies the following elliptic equation:
$${\rm}(\mu(\rho)D\p_{t} v)+\n(\lambda(\rho){\rm div}\p_{t}v)=\n \p_{t}P(\rho)-{\rm}(\p_{t}\mu(\rho)D v)+\n(\p_{t}\lambda(\rho){\rm div}v).$$
We can in a similar way get estimates on $\p_{t}v$ in function of $q$ and $u$.
The rest of the proof is exactly similar to the proof of theorem \ref{theo1} and is nothing than tedious verifications. It is left to the reader.
%_{\widetilde{L}^{1}(B^{\N-1,}_{2,1})}\leq \|T_{q^{n}}({\rm div}u^{n})\|_{\widetilde{L}^{1}(\widetilde{B}^{\N-1,\NN-1}_{2,p,1})}+
%\|T_{{\rm div}u^{n}}q^{n}\|_{\widetilde{L}^{1}(\widetilde{B}^{\N-1,\NN}_{2,p,1})}\\
%&\hspace{8cm}+\|\big(R(q^{n},{\rm div}u^{n})\big)_{BF}\|_{\widetilde{L}^{1}(B^{\N-1}_{2,1})},\\
%&\leq C\big(\|q^{n}\|_{\widetilde{L}^{2}(\widetilde{B}^{\N,\NN}_{2,p,1})}
%\|u^{n}\|_{\widetilde{L}^{2}(\widetilde{B}^{\N,\frac{N}{p}}_{2,p,1})}+\|q^{n}\|_{\widetilde{L}^{\infty}(\widetilde{B}^{\N-1,\NN}_{2,p,%1})}
%\|u^{n}\|_{\widetilde{L}^{1}(\widetilde{B}^{\N+1,\frac{N}{p}+1}_{2,p,1})}\big).
\section{Appendix}
This section is devoted to the proof of proposition \ref{hybrid} and of commutators estimates which have been used in section $2$ and $3$. They are based on
paradifferentiel calculus, a tool introduced by J.-M. Bony in \cite{BJM}. The basic idea of paradifferential calculus is that
any product of two distributions $u$ and $v$ can be formally decomposed into:
$$uv=T_{u}v+T_{v}u+R(u,v)=T_{u}v+T^{'}_{v}u$$
where the paraproduct operator is defined by $T_{u}v=\sum_{q}S_{q-1}u\D_{q}v$, the remainder operator $R$, by
$R(u,v)=\sum_{q}\D_{q}u(\D_{q-1}v+\D_{q}v+\D_{q+1}v)$ and $T^{'}_{v}u=T_{v}u+R(u,v)$.
\begin{proposition}
\label{hybrid}
Let $p_{1},p_{2},p_{3},p_{4}\in[1,+\infty], (s_{1},s_{2},s_{3},s_{4})\in\R^{4}$ and $(p,q)\in[1,+\infty]^{2}$, we have then the following inequalities:
\begin{itemize}
\item If $\frac{1}{p}\leq\frac{1}{p_{2}}+\frac{1}{\lambda}\leq 1$, $\frac{1}{q}\leq\frac{1}{p_{4}}+\frac{1}{\lambda^{'}}\leq 1$ with $(\lambda,\lambda^{'})\in[1,+\infty]^{2}$ and $p_{1}\leq\lambda^{'}$, $p_{1}\leq\lambda$, $p_{3}\leq\lambda^{'}$ then:
\begin{equation}
\|T_{u}v\|_{\widetilde{B}^{s_{1}+s_{2}+\NN-\frac{N}{p_{1}}-\frac{N}{p_{2}}
,s_{3}+s_{4}+\frac{N}{q}-\frac{N}{p_{3}}-\frac{N}{p_{4}}}_{p,q,r}}\lesssim \|u\|_{\widetilde{B}^{s_{1},s_{3}}_{p_{1},p_{3},1}}\|v\|_{\widetilde{B}^{s_{2},s_{4}}_{p_{2},p_{4},r}},
\label{62}
\end{equation}
if $s_{1}+\frac{N}{\lambda^{'}}\leq\frac{N}{p_{1}}$, $s_{1}+\frac{N}{\lambda}\leq\frac{N}{p_{1}}$ and $s_{3}+\frac{N}{\lambda^{'}}\leq\frac{N}{p_{3}}$.
\item If %$\frac{1}{p}\leq\frac{1}{p_{1}}+\frac{1}{p_{2}}$,
$\frac{1}{q}\leq\frac{1}{p_{3}}+\frac{1}{p_{4}}$ and
%$s_{1}+s_{2}+N\inf(0,1-\frac{1}{p_{1}}-\frac{1}{p_{2}})>0,
$s_{3}+s_{4}+N\inf(0,1-\frac{1}{p_{3}}-\frac{1}{p_{4}})>0$
then
 \begin{equation}
 \sum_{l\geq 4}2^{l(s_{3}+s_{4}+\frac{N}{q}-\frac{N}{p_{3}}-\frac{N}{p_{4}})}\|\D_{l}R(u,v)\|_{L^{q}}\lesssim \|u\|_{\widetilde{B}^{s_{1},s_{3}}_{p_{1},p_{3},1}}\|v\|_{\widetilde{B}^{s_{2},s_{4}}_{p_{2},p_{4},r}}.
\label{63}
\end{equation}
\item If %$\frac{1}{p}\leq\frac{1}{p_{1}}+\frac{1}{p_{2}}$,
$\frac{1}{p}\leq\frac{1}{p_{3}}+\frac{1}{p_{4}}\leq 1$, $\frac{1}{p}\leq\frac{1}{p_{3}}+\frac{1}{p_{2}}\leq 1$, $\frac{1}{p}\leq\frac{1}{p_{1}}+\frac{1}{p_{4}}\leq 1$, $\frac{1}{p}\leq\frac{1}{p_{1}}+\frac{1}{p_{2}}\leq 1$ and
%$s_{1}+s_{2}+N\inf(0,1-\frac{1}{p_{1}}-\frac{1}{p_{2}})>0,
$s_{3}+s_{4}>0$, $s_{3}+s_{2}>0$, $s_{4}+s_{1}>0$, $s_{1}+s_{2}>0$
then
 \begin{equation}
 \sum_{l\leq 4}2^{l(s_{1}+s_{2}+\frac{N}{p}-\frac{N}{p_{1}}-\frac{N}{p_{2}})}\|\D_{l}R(u,v)\|_{L^{p}}\lesssim \|u\|_{\widetilde{B}^{s_{1},\frac{N}{p_{3}}-\frac{N}{p_{1}}+s_{1}}_{p_{1},p_{3},1}}\|v\|_{\widetilde{B}^{s_{2},\frac{N}{p_{4}}-\frac{N}{p_{2}}+s_{2}}_{p_{2},p_{4},r}}.
\label{63}
\end{equation}
with $s_{3}=\frac{N}{p_{3}}-\frac{N}{p_{1}}+s_{1}$ and $s_{4}=\frac{N}{p_{4}}-\frac{N}{p_{2}}+s_{2}$.
%\item If $\frac{1}{p}\leq\frac{1}{p_{1}}+\frac{1}{p_{2}}\leq 1$, $\frac{1}{q}\leq\frac{1}{p_{3}}+\frac{1}{p_{4}}\leq1$ and
%$s_{1}+s_{2}=s_{3}+s_{4}=0$.
% then:
% \begin{equation}
% \|R(u,v)\|_{\widetilde{B}^{\NN-\frac{N}{p_{1}}-\frac{N}{p_{2}}
%,\frac{N}{q}-\frac{N}{p_{3}}-\frac{N}{p_{4}}}_{p,q,\infty}}\lesssim \|u\|_{\widetilde{B}^{s_{1},s_{3}}_{p_{1},p_{3},1}}\|v\|_{\widetilde{B}^{s_{2},s_{4}}_{p_{2},p_{4},\infty}}.
%\label{64}
%\end{equation}
\item If $u\in L^{\infty}$, we also have:
\begin{equation}
\|T_{u}v\|_{\widetilde{B}^{s_{1}
,s_{2}}_{p,q,r}}\lesssim \|u\|_{L^{\infty}}\|v\|_{\widetilde{B}^{s_{1},s_{2}}_{p,q,r}},
\label{65}
\end{equation}
and if $\min(s_{1},s_{2})>0$ then:
\begin{equation}
\|R(u,v)\|_{\widetilde{B}^{s_{1}
,s_{2}}_{p,q,r}}\lesssim \|u\|_{L^{\infty}}\|v\|_{\widetilde{B}^{s_{1},s_{2}}_{p,q,r}}.
\label{66}
\end{equation}
\end{itemize}
\end{proposition}
{\bf Proof:} Let us prove (\ref{62}). According to the decomposition of J.-M. Bony \cite{BJM}, we have:
$$uv=T_{u}v+T_{v}u+R(u,v),$$
so for all $l>0$:
$$\D_{l}T_{u}v=\sum_{|l-l^{'}|\leq3}\D_{l}(S_{l^{'}-1}u\D_{l^{'}}v),$$
For $\alpha,\beta\in\R$, let us define the following characteristic function on $\mathbb{Z}$
$$
\begin{aligned}
\va^{\alpha,\beta}=\alpha\;\;\;\mbox{if}\;\;r\leq0,\\
\va^{\alpha,\beta}=\beta\;\;\;\mbox{if}\;\;r\geq1.
\end{aligned}
$$
if  $\frac{1}{p}\leq\frac{1}{p_{2}}+\frac{1}{\lambda}\leq 1$ and $\frac{1}{q}\leq\frac{1}{p_{4}}+\frac{1}{\lambda^{'}}\leq 1$ then
$$\|\D_{l}T_{u}v\|_{L^{\va^{p,q}(l)}}\lesssim 2^{lN\va^{\frac{1}{p_{2}}+\frac{1}{\lambda}-\frac{1}{p},\frac{1}{p_{4}}+\frac{1}{\lambda^{'}}-\frac{1}{q}}(l)}
\sum_{|l-l^{'}|\leq3}\|S_{l^{'}-1}u\|_{L^{\va^{\lambda,\lambda^{'}}(l^{'})}}\|
\D_{l^{'}}v\|_{L^{\va^{p_{2},p_{4}}(l^{'})}}.$$
We have by Berstein inequalities and as $p_{1}\leq\lambda^{'}$, $p_{3}\leq\lambda^{'}$, $p_{1}\leq\lambda$ and
$s_{1}+\frac{N}{\lambda}\leq \frac{N}{p_{1}}$, $s_{1}+\frac{N}{\lambda^{'}}\leq \frac{N}{p_{1}}$, $s_{3}+\frac{N}{\lambda^{'}}\leq \frac{N}{p_{3}}$:
$$
\begin{aligned}
\|S_{l^{'}-1}u\|_{L^{\va^{\lambda,\lambda^{'}}(l^{'})}}&\lesssim\sum_{k\leq l^{'}-2}2^{k(\va^{\frac{N}{p_{1}},\frac{N}{p_{3}}}(k)-\va^{\frac{N}{\lambda},\frac{N}{\lambda^{'}}}(l^{'}))}
\|\D_{k}u\|_{L^{\va^{p_{1},p_{3}}(k)}}\\
&\lesssim\sum_{k\leq l^{'}-2}2^{k(\va^{\frac{N}{p_{1}}-s_{1},\frac{N}{p_{3}}-s_{3}}(k)-\va^{\frac{N}{\lambda},\frac{N}{\lambda^{'}}}(l^{'}))}2^{k\va^{s_{1},s_{3}}(k)}
\|\D_{k}u\|_{L^{\va^{p_{1},p_{3}}(k)}}\\
&\lesssim 2^{l^{'}(\va^{\frac{N}{p_{1}}-s_{1},\frac{N}{p_{3}}-s_{3}}(l^{'})-\va^{\frac{N}{\lambda},\frac{N}{\lambda^{'}}}(l^{'}))}
\|u\|_{\widetilde{B}^{s_{1},s_{3}}_{p_{1},p_{3},1}}.
\end{aligned}
$$
Since $\|\D_{l^{'}}v\|_{L^{\va^{p_{2},p_{4}}(l^{'})}}=c_{l^{'}}2^{-l^{'}(\va^{s_{2},s_{4}}(l^{'}))}
\|v\|_{\widetilde{B}^{s_{2},s_{4}}_{p_{2},p_{4},1}}$ with $\sum_{l^{'}\in\mathbb{Z}}c_{l^{'}}\leq1$
we finally gather as $l>0$:
$$
\begin{aligned}
\|\D_{l}T_{u}v\|_{L^{q}}\lesssim c_{l}2^{l \va^{\frac{N}{p_{1}}+\frac{N}{p_{2}}-\frac{N}{p}-s_{1}-s_{2},
\frac{N}{p_{2}}+\frac{N}{p_{4}}-\frac{N}{q}-s_{3}-s_{4}}(l)}
\|u\|_{\widetilde{B}^{s_{1},s_{3}}_{p_{1},p_{3},1}}\|v\|_{\widetilde{B}^{s_{2},s_{4}}_{p_{2},p_{4},1}}.
%\lesssim c_{l}2^{l\big(\frac{N}{p_{4}}+\frac{N}{p_{3}}-\frac{N}{q}-s_{3}-s_{4}\big)}
%\|u\|_{\widetilde{B}^{s_{1},s_{3}}_{p_{1},p_{3},1}}\|v\|_{\widetilde{B}^{s_{2},s_{4}}_{p_{2},p_{4},1}}.
%\label{111}
\end{aligned}
$$
And we obtain (\ref{62}).
\\
\\
Straightforward modification give (\ref{65}). In this case as $\|S_{k-1}u\|_{L^{\infty}}\leq \|u\|_{L^{\infty}}$ we have:
$$\|\D_{l}T_{u}v\|_{L^{\va^{p,q}(l)}}\lesssim \sum_{|l-l^{'}|\leq3}\|u\|_{L^{\infty}}\|
\D_{l^{'}}v\|_{L^{\va^{p_{2},p_{4}}(l^{'})}}.$$
Next we have:
$$2^{l\va^{p_{2},p_{4}}(l)}\|\D_{l}T_{u}v\|_{L^{\va^{p,q}(l)}}\lesssim \|u\|_{L^{\infty}} \sum_{|l-l^{'}|\leq3}2^{l\va^{p_{2},p_{4}}(l)-l^{'}\va^{p_{2},p_{4}}(l^{'}))}2^{\va^{p_{2},p_{4}}(
l^{'})}\|
\D_{l^{'}}v\|_{L^{\va^{p_{2},p_{4}}(l^{'})}}.$$
We conclude by convolution.\\
\\
To prove (\ref{63}), we write:
$$\D_{l}R(u,v)=\sum_{k\geq l-2}\D_{l}(\D_{k}u\widetilde{\D}_{k}v).$$
We consider now the case $l>3$. By Bernstein and H\"older inequalities we obtain when  $\frac{1}{q}\leq \frac{1}{p_{3}}+\frac{1}{p_{4}}\leq 1$:
$$\|\D_{l}R(u,v)\|_{L^{q}}\lesssim2^{Nl(\frac{1}{p_{3}}+\frac{1}{p_{4}}-\frac{1}{q})}\sum_{k\geq l-2}
\|\D_{k}u\|_{L^{p_{3}}}\|\widetilde{\D}_{k}v\|_{L^{p_{4}}}.$$
Next we have:
$$
\begin{aligned}
2^{l(s_{3}+s_{4}+\frac{N}{q}-\frac{N}{p_{3}}-\frac{N}{p_{4}})}\|\D_{l}R(u,v)\|_{L^{q}}&\lesssim
\sum_{k\geq l-2}
2^{(l-k)(s_{3}+s_{4})}2^{ks_{3}}\|\D_{k}u\|_{L^{p_{3}}}2^{ks_{4}}\|\widetilde{\D}_{k}v\|_{L^{p_{4}}},\\
&\lesssim (c_{k})*(d_{k^{'}}),
\end{aligned}
$$
with $c_{k}=1_{[-\infty,2]}(k)2^{k(s_{3}+s_{4})}$ and $d_{k^{'}}=2^{k^{'}s_{3}}\|\D_{k}u\|_{L^{p_{3}}}2^{k^{'}s_{4}}\|\widetilde{\D}_{k}v\|_{L^{p_{4}}}$
. We conclude by Young inequality as $s_{3}+s_{4}>0$.\\
\\
We have to treat now the case when $l<0$. We have then as $\frac{1}{p}\leq\frac{1}{p_{3}}+\frac{1}{p_{4}}\leq 1$ and
$\frac{1}{p}\leq\frac{1}{p_{1}}+\frac{1}{p_{2}}\leq 1$:
$$
\begin{aligned}
&\|\D_{l}R(u,v)\|_{L^{p}}\lesssim2^{Nl(\frac{1}{p_{3}}+\frac{1}{p_{4}}-\frac{1}{p})}\sum_{k\geq 2}
\|\D_{k}u\|_{L^{p_{3}}}\|\widetilde{\D}_{k}v\|_{L^{p_{4}}}\\
&\sum_{0\leq k\leq 1, |k-k^{'}|\leq 1}
\|\D_{k}u\D_{k^{'}}v\|_{L^{p}}+2^{Nl(\frac{1}{p_{1}}
+\frac{1}{p_{2}}-\frac{1}{p})}\sum_{l-2\leq k\leq -1}
\|\D_{k}u\|_{L^{p_{1}}}\|\widetilde{\D}_{k}v\|_{L^{p_{2}}}.
\end{aligned}
$$
And by convolution on the middle frequencies:
$$
\begin{aligned}
&\|\D_{l}R(u,v)\|_{L^{p}}\lesssim2^{Nl(\frac{1}{p_{3}}+\frac{1}{p_{4}}-\frac{1}{p})}\sum_{k\geq 2}
\|\D_{k}u\|_{L^{p_{3}}}\|\widetilde{\D}_{k}v\|_{L^{p_{4}}}\\
&(2^{l(\frac{N}{p_{3}}
+\frac{N}{p_{2}}-\frac{N}{p}-s_{3}-s_{2})}+2^{l(\frac{N}{p_{1}}
+\frac{N}{p_{4}}-\frac{N}{p}-s_{1}-s_{4})}
)c_{l}+
2^{Nl(\frac{1}{p_{1}}
+\frac{1}{p_{2}}-\frac{1}{p})}\sum_{l-2\leq k\leq -1}
\|\D_{k}u\|_{L^{p_{1}}}\|\widetilde{\D}_{k}v\|_{L^{p_{2}}},
\end{aligned}
$$
with $c_{l}\in l^{1}(\mathbb{Z})$.
%We recall now that: $p_{1}\geq p_{2}$ and $p_{3}\geq p_{4}$, we have then by Bernstein inequalities:
Next by convolution we obtain:
$$
\begin{aligned}
&\|\D_{l}R(u,v)\|_{L^{p}}\lesssim c_{l}(2^{l(\frac{N}{p_{3}}+\frac{N}{p_{4}}-\frac{N}{p}-s_{3}-s_{4})}+2^{l(\frac{N}{p_{3}}
+\frac{N}{p_{2}}-\frac{N}{p}-s_{3}-s_{2})}+2^{l(\frac{N}{p_{1}}
+\frac{N}{p_{4}}-\frac{N}{p}-s_{1}-s_{4})}\\
&\hspace{6cm}+
2^{l(\frac{N}{p_{1}}+\frac{N}{p_{2}}-\frac{N}{p}-s_{1}-s_{2})})
\|u\|_{\widetilde{B}^{s_{1},s_{3}}_{p_{1},p_{3},1}}\|v\|_{\widetilde{B}^{s_{2},s_{4}}_{p_{2},p_{4},r}}.
\end{aligned}
$$
And we can conclude.
\\
We want prove now the inequality (\ref{65}). We have then:
$$
\begin{aligned}
2^{l\va^{s_{1},s_{2}}(l)}\|\D_{l}R(u,v)\|_{L^{p}}&\lesssim\sum_{k\geq l-2}2^{(l-k)\va^{s_{1},s_{2}}(l)}
2^{k\va^{s_{1},s_{2}}(l)}\|\D_{k}u\|_{L^{\infty}}\|\widetilde{\D}_{k}v\|_{L^{\va^{p,q}(k)}},
\end{aligned}
$$
And we conclude by Young inequality.
{\hfill $\Box$}\\
\\
%\begin{lemme}
%\label{alemme1}
%Let a be a bounded continuous function such that $a\geq\bar{a}>0$. Let $N\geq 1$, $1<p<+\infty$ and $u\in L^{p}(\R^{N})$
% such that ${\rm supp}\hat{u}\subset C(0,R_{1},R_{2})$ (with $0<R_{1}<R_{2}$). Then we have
%for a constant $c$ depending only on $N$ and $\frac{R_{2}}{R_{1}}$:
%$$c\bar{a}\frac{p-1}{p^{2}}R_{1}^{2}\int_{\R^{N}}|u|^{p}dx\leq -\int_{\R^{N}}{\rm div}(a\n u)|u|^{p-2}udx.$$
%\end{lemme}
%Inequalities (\ref{12}) and (\ref{18}) are consequence of the following lemma:
\begin{lemme}
\label{alemme2}
Let $1\leq p_{1}\leq p\leq+\infty$ and $\sigma\in(-\min(\NN,\frac{N}{p_{1}^{'}}),\NN+1]$. There exists a sequence $c_{q}\in l^{1}(\mathbb{Z})$ such that $\|c_{q}\|_{l^{1}}=1$
and a constant
$C$ depending only on $N$ and $\sigma$ such that:
\begin{equation}
\forall q\in\mathbb{Z},\;\;\|[v\cdot\n,\D_{q}]a\|_{L^{p_{1}}}\leq C c_{q}2^{-q\sigma}\|\n v\|_{B^{\NN}_{p,1}}
\|a\|_{B^{\sigma}_{p_{1},1}}.
\label{52}
\end{equation}
In the limit case $\sigma=-\min(\NN,\frac{N}{p_{1}^{'}})$, we have:
\begin{equation}
\forall q\in\mathbb{Z},\;\;\|[v\cdot\n,\D_{q}]a\|_{L^{p_{1}}}\leq C c_{q}2^{q\NN}\|\n v\|_{B^{\NN}_{p,1}}
\|a\|_{B^{-\frac{N}{p_{1}}}_{p,\infty}}.
\label{53}
\end{equation}
Finally, for all $\sigma>0$ and $\frac{1}{p_{2}}=\frac{1}{p_{1}}-\frac{1}{p}$, there exists a constant $C$ depending only on $N$ and on $\sigma$ and a sequence
$c_{q}\in l^{1}(\mathbb{Z})$ with norm
$1$ such that:
\begin{equation}
\forall q\in\mathbb{Z},\;\;\|[v\cdot\n,\D_{q}]v\|_{L^{p}}\leq C c_{q}2^{-q\sigma}(\|\n v\|_{L^{\infty}}\|v\|_{B^{\sigma}_{p_{1},1}}+\|\n v\|_{L^{p_{2}}}\|\n v\|_{B^{\sigma-1}_{p,1}}).
\label{54}
\end{equation}
\end{lemme}
{\bf Proof:} These results are proved in \cite{BCD} chapter $2$.
%Inequality (\ref{52}) has been proved in (\cite{DD}), lemma A1 under the hypothesis (which does not play any role
%if $\sigma>\NN$) that ${\rm div}v=0$. It is based on the
%decomposition:
%\begin{equation}
%[v\cdot\n,\D_{q}]a=[T_{v_{j}},\D_{q}]\p_{j}a+T^{'}_{\D_{q}\p_{j}a}v^{j}-
%\D_{q}T_{\p_{j}a}v^{j}-\D_{q}R(\p_{j}a,v^{j}).
%\label{55}
%\end{equation}
%In the case $\sigma=-\NN$, the computations made in \cite{DD} show that the $L^{p}$ norm of the first three terms in (\ref{55})
%may be bounded
%by $C2^{q\NN}\|\n v\|_{B^{\NN}_{p,1}}\|a\|_{B^{-\NN}_{p,\infty}}$.\\
%For bounding the last term, we use the following classical result of continuity for the remainder (see \cite{RS}):
%\begin{equation}
%\|R(f,g)\|_{B^{-\N}_{2,\infty}}\lesssim\|f\|_{B^{-s}_{2,\infty}}\|g\|_{B^{s}_{2,1}}
%\label{56}
%\end{equation}
%which holds for all real number $s$. This yields (\ref{53}).\\
%The proof of (\ref{54}) relies on a similar arguments. The details are left to the reader.
{\hfill $\Box$}\\
%Inequality (\ref{13}) is a consequence of the following lemma:

%Appliquer la chaleur sur les solutions \`a
%Bresch,Desjardins avec une bonne hypoth\`ese sur la pression,
%on devrait obtenir de l'unicit\'e au moins en temps petit,
%voir pour des temps plus grand si on peut avoir des methode type Gronwall.

\begin{thebibliography}{}
\bibitem{AP}
H. Abidi and M. Paicu. \'Equation de Navier-Stokes avec densit\'e et viscosit\'e variables dans l'espace critique. \textit{Annales de l'institut Fourier}, 57 no. 3 (2007), p. 883-917.
\bibitem{BC}
H. Bahouri and J.-Y. Chemin, \'Equations d'ondes quasilin\'eaires et
estimation de Strichartz, \textit{Amer. J. Mathematics}. 121 (1999) 1337-1377.
\bibitem{BCD}
H. Bahouri, J.-Y. Chemin and R. Danchin. Fourier analysis and nonlinear partial differential equations,
\textit{to appear in Springer}.
\bibitem{BJM}
J.-M. Bony, Calcul symbolique et propagation des singularit\'es pour
les \'equations aux d\'eriv\'ees partielles non lin\'eaires, \textit{Annales
Scientifiques de l'\'ecole Normale Sup\'erieure}. 14 (1981)
209-246.
\bibitem{BD}
D. Bresch and B. Desjardins, Existence of global weak solutions for
a 2D Viscous shallow water equations and convergence to the
quasi-geostrophic model. \textit{Comm. Math. Phys.}, 238(1-2): 211-223,
2003.
\bibitem{CD}
F. Charve and R. Danchin, Global existence in critical spaces for compressible Navier-Stokes equations, preprint and submitted.
\bibitem{CT}
J.-Y. Chemin, Th\'eor\`emes d'unicit\'e pour le syst\`eme de
Navier-Stokes tridimensionnel, \textit{J.d'Analyse Math}. 77 (1999) 27-50.
\bibitem{CA}
J.-Y. Chemin, About Navier-Stokes system, \textit{Pr\'epublication du
Laboratoire d'Analyse Num\'erique de Paris 6}, R96023 (1996).
\bibitem{CL}
J.-Y. Chemin and N. Lerner, Flot de champs de vecteurs non
lipschitziens et \'equations de Navier-Stokes, \textit{J.Differential
Equations}, 121 (1992) 314-328.
\bibitem{CMZ}
Q. Chen, C. Miao and Z. Zhang,
Global well-posedness for the compressible Navier-Stokes equations with the highly oscillating initial velocity, \textit{arXiv:0907.4540v2}.
\bibitem{DFourier}
R. Danchin, Fourier analysis method for PDE's, \textit{Preprint}, Novembre 2005.
\bibitem{DL}
R. Danchin, Local Theory in critical Spaces for Compressible Viscous
and Heat-Conductive Gases, \textit{Communication in
Partial Differential Equations}, 26 (78),1183-1233, (2001).
\bibitem{DG}
R. Danchin, Global existence in critical spaces for compressible Navier-Stokes equations, \textit{Inventiones Mathematicae}, 141, pages 579-614 (2000).
\bibitem{DG1}
R. Danchin, Global existence in critical spaces for compressible viscous and heat-conductive gases, \textit{Archiv for Rational Mechanics and Analysis}, 160, pages 1-39 (2001).
%\bibitem{DG}
%R. Danchin, Global Existence in Critical Spaces for Flows of
%Compressible Viscous and Heat-Conductive Gases, \textit{Arch.Rational
%Mech.Anal}.160, (2001), 1-39.
\bibitem{DU}
R. Danchin, On the uniqueness in critical spaces for compressible Navier-Stokes equations. \textit{NoDEA Nonlinear Differentiel Equations
Appl}, 12(1):111-128, 2005.
\bibitem{DW}
R. Danchin, Well-Posedness in Critical Spaces for Barotropic Viscous Fluids with Truly Not Constant
Density, \textit{Communications in Partial Differential Equations},32:9,1373-1397.
\bibitem{H1}
B. Haspot,Cauchy problem for viscous shallow water equations  with a term of capillarity , \textit{accepted in HYP 2008}.
\bibitem{H}
B. Haspot, Local well-posedness results for density-dependent incompressible fluids, \textit{Arxiv}, 0902.1982 (February 2009).
\bibitem{H2}
B. Haspot, Well-posedness in critical spaces for barotropic viscous fluids, \textit{Arxiv}, (March 2009).
\bibitem{5H1}
D. Hoff. Global existence for 1D, compressible, isentropic
Navier-Stokes equations with large initial data.
\textit{Trans. Amer. Math. Soc}, 303(1), 169-181, 1987.
\bibitem{5H5}
D. Hoff, Uniqueness of weak solutions of the Navier–Stokes equations of multidimensional, compressible flow, \textit{SIAM J. Math. Anal}. 37 (6) (2006).
\bibitem{5H4}
D. Hoff. Discontinuous solutions of the Navier-Stokes equations
for multidimensional flows of the heat conducting fluids.
\textit{Arch. Rational Mech. Anal.}, 139, (1997), p. 303-354.
\bibitem{5H2}
D. Hoff. Global solutions of the Navier-Stokes equations for
multidimensional compressible flow with discontinuous initial
data. \textit{J. Differential Equations}, 120(1), 215-254, 1995.
\bibitem{5H3}
D. Hoff. Strong convergence to global solutions for
multidimensional flows of compressible, viscous fluids with
polytropic equations of state and discontinuous initial data.
\textit{Arch. Rational Mech. Anal.}, 132(1), 1-14, 1995.
%data. J. Differential Equations, 120(1):215-254, 1995.\\
\bibitem{5HZ}
D. Hoff and K. Zumbrum. Multi-dimensional diffusion waves for
the Navier-Stokes equations of compressible flow,
\textit{Indiana University Mathematics Journal}, 1995, 44, 603-676.
\bibitem{5K1}
A. V. Kazhikov. The equation of potential flows of a compressible
viscous fluid for small Reynolds numbers: existence,
uniqueness and stabilization of solutions. \textit{Sibirsk. Mat. Zh.}, 34 (1993), no. 3, p. 70-80.
\bibitem{5K}
A. V. Kazhikov and V. V. Shelukhin. Unique global solution with
respect to time of initial-boundary value problems for one-
dimensional equations of a viscous gas. \textit{Prikl. Mat. Meh.}, 41(2): 282-291, 1977.
\bibitem{MN1}
Akitaka Matsumura and Takaaki Nishida. The initial value problem for
the equations of motion of compressible
viscous and heat-conductive gases. \textit{J. Math. Kyoto Univ.}, 20(1): 67-104, 1980.
\bibitem{MN}
Akitaka Matsumura and Takaaki Nishida. The initial value problem for
the equations of motion of compressible
viscous and heat-conductive fluids. \textit{Proc. Japan Acad. Ser. A Math. Sci}, 55(9):337-342, 1979.
\bibitem{Meyer}
Y. Meyer. Wavelets,paraproducts, and Navier-Stokes equation. \textit{In Current developments in mathematics, 1996 (Cambridge, MA)},
page 105-212. Int. Press, Boston, MA, 1997.
\bibitem{Nash}
J. Nash. Le probl\`eme de Cauchy pour les \'equations diff\'erentielles d'un fluide g\'en\'eral. \textit{Bull. Soc. Math. France}, 90:
487-497, 1962.
\bibitem{RS}
T. Runst and W. Sickel: Sobolev spaces of fractional order, Nemytskij operators, and nonlinear partial differential equations. \textit{de Gruyter Series in
Nonlinear Analysis and Applications}, 3. Walter de Gruyter and Co., Berlin (1996)\\
na  21, no. 1 (2005), 1-24.
\bibitem{5S}
D. Serre. Solutions faibles globales des \'equations de
Navier-Stokes pour un fluide compressible.\textit{Comptes rendus de l'Acad\'emie des sciences}. S\'erie 1,
303(13): 639-642, 1986.
\bibitem{5So}
V. A. Solonnikov. Estimates for solutions of nonstationary
Navier-Stokes systems. \textit{Zap. Nauchn. Sem. LOMI}, 38,
(1973), p.153-231; J. Soviet Math. 8, (1977), p. 467-529.
\bibitem{5V}
V. Valli and  W. Zajaczkowski. Navier-Stokes equations for compressible
fluids: global existence and qualitative properties of the solutions
in the general case. \textit{Commun. Math. Phys.}, 103, (1986) no 2, p.
259-296.
\end{thebibliography}
\end{document}